\documentclass[a4paper,12pt]{amsart}

\usepackage{amssymb}
\usepackage{latexsym} 
\usepackage{amsfonts} 
\usepackage{amsmath}
\usepackage{eucal} 
\usepackage{bm} 
\usepackage{bbm} 
\usepackage{graphicx} 
\usepackage[english]{varioref} 
\usepackage[nice]{nicefrac} 
\usepackage[all]{xy}
\usepackage{amsthm}

\def\ge{\geqslant}
\def\le{\leqslant}
\def\a{\alpha}
\def\b{\beta}

\def\d{\delta}

\def\e{\epsilon}

\def\o{\omega}

\def\s{\sigma}
\def\t{\tau}
\def\th{\theta}

\def\l{\lambda}

\def\i{^{-1}}
\def\ZZ{\mathbb Z}
\def\NN{\mathbb N}
\def\QQ{\mathbb Q}
\def\RR{\mathbb R}

\def\FF{\mathbb F}
\newcommand{\kk}{\Bbbk}

\def\cb{\mathcal B}

\def\cd{\mathcal D}

\def\co{\mathcal O}
\def\cp{\mathcal P}

\def\aa{\mathbf a}

\def\tz{\tilde z}
\def\tl{\tilde \l}

\def\ty{\tilde y}

\def\tW{\tilde W}

\def\tS{\tilde S}

\def\aa{\mathbf a}

\theoremstyle{plain}
\newtheorem{thm}{Theorem}[section] 
\newtheorem*{thm*}{Theorem} 
 \newtheorem{prop}[thm]{Proposition}
 \newtheorem{lem}[thm]{Lemma}
 \newtheorem{cor}[thm]{Corollary}

\theoremstyle{definition}

\theoremstyle{remark}

\newtheorem*{rmk}{Remark}
\newtheorem*{claim*}{Claim}

\begin{document}

\author{Xuhua He}
\address{Department of Mathematics, The Hong Kong University of Science and Technology, Clear Water Bay, Kowloon, Hong Kong}
\thanks{The author is partially supported by HKRGC grants
601409.}
\email{maxhhe@ust.hk}
\title[]{On some partitions of an affine flag variety}

\begin{abstract}
In this paper, we discuss some partitions of affine flag varieties. These partitions include as special cases the partition of affine flag variety into affine Deligne-Lusztig varieties and the affine analogue of the partition of flag varieties into $\cb_w(b)$ introduced by Lusztig in \cite{L1} as part of the definition of character sheaves. 

Among other things, we give a formula for the dimension of affine Deligne-Lusztig varieties for classical loop groups in terms of degrees of class polynomials of extended affine Hecke algebra. We also prove that any simple $GL_n(\FF_q((\e)))$-module occurs as a subquotient of the cohomology of affine Deligne-Lusztig variety $X_w(1)$ for some $w$ in the extended affine Weyl group $\ZZ^n \rtimes S_n$ must occurs for some $w$ in the finite Weyl group $S_n$. Similar result holds for $Sp_{2n}$. 
\end{abstract}

\maketitle

\section*{Introduction}

\subsection{} Let $G$ be a connected reductive algebraic group over an algebraically closed field $\kk$. We consider a ``twisted'' conjugation action on $G$ defined by $g \cdot_\s g'=g g' \s(g) \i$. Here $\s$ is the identity map or a Frobenius morphism on $G$ (in case where $\kk$ is of positive characteristic). 

Let $B$ be a Borel subgroup of $G$ with $\s(B)=B$ and $W$ be the Weyl group of $G$. For $b \in G$ and $w \in W$, we set $$\cb_{w, \s}(b)=\{g B \in G/B; g \i b \s(g) \in B \dot w B\}.$$ Then we have a partition of flag variety $G/B=\sqcup_{w \in W} \cb_{w, \s}(b)$. 

In the case where $\s$ is a Frobenius morphism, $\cb_{w, \s}$ is a Deligne-Lusztig variety introduced in \cite{DL}. In the case where $\s$ is the identity, $\cb_{w, \s}(b)$ appears as part of the definition of character sheaves introduced in \cite{L1}. These varieties and their cohomology groups play an important role in geometric representation theory. 

\subsection{} Let $L=\kk((\e))$ is the field of formal Laurent power series and $o=\kk[[\e]]$ be the ring of formal power series. Let $I$ be a Iwahori subgroup of the loop group $G(L)$. The quotient $G(L)/I$ is called an affine flag variety. For some technical reason, we choose $I$ to be the inverse image of $B^-$ under the projection map $G(o) \to G$ sending $\e$ to $0$, here $B^-$ is a Borel subgroup of $G$ opposite to $B$. 

The main purpose of this paper is to study an analogue of the above partition in the affine case and their cohomology groups. 

We consider a ``twisted'' conjugation action of $G(L)$ on itself as $g \cdot_\s h=g h \s(g) \i$ for $g, h \in G(L)$. Here $\s$ is a bijective group homomorphism on $G(L)$ of one of the following type:

(1) For any nonzero element $a \in \kk$, define $\s_a(p(t))=p(a \cdot t)$ for any formal Laurent power series $p(t)$. We extend $\s_a$ to a group homomorphism on $G(L)$, which we still denote by $\s_a$.

(2) If $\kk$ is of positive characteristic and $F: \kk \to \kk$ be a Frobenius morphism. Then set $F(\sum a_n t^n)=\sum F(a_n) t^n$. We extend $F$ to a group homomorphism on $G(L)$, which we denote by $\s_F$. 

The $\s_a$-conjugacy classes are studies by Baranovsky and Ginzburg in \cite{BG96}. The $\s_F$-conjugacy classes are studied by Kottwitz in \cite{Ko97}. 

\subsection{} Let $\tW$ be the extended affine Weyl group of $G(L)$. For $b \in G(L)$ and $w \in \tW$, we set $$X_{w, \s}(b)=\{g I \in G(L)/I; g \i b \s(g) \in I \dot w I\}.$$ Then we have a partition of affine flag variety $$G(L)/I=\sqcup_{w \in \tW} X_{w, \s}(b).$$ In the case that $\s=\s_F$, $X_{w, \s}(b)$ is called an {\it affine Deligne-Lusztig variety}. In the case that $\s$ is identity map, $X_{1, \s}(b)$ is called an {\it affine Springer fiber}. $X_{w, id}(b)$ is also considered in \cite[Section 7]{L3}. 

A big difference between the finite case discussed in subsection 0.1 and its affine analogue above is that $X_{w, \s}(b)$ is not always nonempty for $\s$ coming from Frobenius morphism. A challenging problem is to determine the empty/nonempty pattern for affine Deligne-Lusztig varieties. For more details, see the discussions in \cite{GHKR} and \cite{GH}. 

Let $Z_{G(L), \s}(b)=\{g \in G(L); g \i b \s(g)=b\}$ be the centralizer of $b$ for the twisted conjugation action. Then $Z_{G(L), \s}(b)$ acts on $X_{w, \s}(b)$ in the natural way. If moreover, $X_{w, \s}(b)$ is finite-dimensional, then $Z_{G(L), \s}(b)$ also acts on $H^*_c(X_{w, \s}(b))$. Any simple $Z_{G(L), \s}(b)$-modules occurs as a subquotient of $H^*_c(X_{w, \s}(b))$ is called {\it obtained from cohomological construction}. 

\subsection{} Our starting point is the following stratification of $G(L)$ into locally closed subschemes that are equivariant for the ``twisted'' conjugation action $$G(L)=\sqcup_{[w] \in \tW_{good}/\approx} G(L) \cdot_\s I \dot w I.$$ In the case that $\s=\s_F$, this follows from Kottwitz's classification of $\s_F$-conjugacy classes and several properties about the ``good'' elements of extended affine Weyl groups established in section 1 (see subsection \ref{reK}). In the case that $\s=\s_a$, the stratification is proved in Prop \ref{str} for $PGL_n(L)$ and the identity component of $PSP_{2n}(L)$. We expect that such stratification holds for any adjoint group $G$. 

We introduce the class polynomials for classical extended affine Hecke algebras, generalization the construction of Geck and Rouquier for finite Hecke algebras in \cite{GR}. Then we prove in Theorem \ref{class} that for classical groups, the nonemptiness (resp. dimension) of affine Deligne-Lusztig varieties corresponds to the nonzeroness (resp. degree) of certain class polynomials. 

For $PGL_n(L)$ and the identity component of $PSP_{2n}(L)$, we obtain a sharper formula in Prop \ref{class3}, $$\dim (X_{w, \s}(b))=\max_{C} \frac{1}{2}(l(w)+l(C)+\deg(f_{w, C}))-l(f \i(\aa)).$$ We also prove that for (i) $\s=\s_a$ with $a$ not a root of unity or (ii) $\s=id$ and $b$ is a regular semisimple integral element, $X_{w, \s}(b)$ is always finite dimensional.

\subsection{} Given $b \in G(L)$, in general, there are infinitely many $w \in \tW$ with $X_{w, \s}(b) \neq \emptyset$. Then a priori, to get all the simple $Z_{G(L), \s}(b)$-modules that are obtained from cohomological construction, one needs to calculate $H^*_c(X_{w, \s}(b))$ for any $w \in \tW$ such that $X_{w, \s}(b) \neq \emptyset$. We will show that for $PGL_n(L)$ and the identity component of $PSP_{2n}(L)$, any simple $Z_{G(L), \s}(b)$-modules that are obtained from cohomological construction must occurs as a subquotient of $H^*_c(X_{w, \s}(b))$, where $w$ is in a finite subset of $\tW$ determined by $b$. The more precise statement can be found in Theorem \ref{main}. 

We mention some applications of this result.

(1) Let $G=PGL_n$. Then any simple $PGL_n(\FF((\e))$-module occurs as a subquotient of the cohomology of affine Deligne-Lusztig variety $X_{w, \s_F}(1)$ for some $w \in \tW$ must occurs for some $w$ in the finite Weyl group $S_n$. See Corollary \ref{finiteA}.

(2) Let $G=PGL_n$ (resp. $G=PSP_{2n}$) and $\chi$ a dominant regular coweight (resp. dominant regular coweight in the coroot lattice). Then any simple $Z_{G(L), \s}(\e^\chi)$-module occurs as a subquotient of $H^*_c(X_{w, \s}(\e^\chi))$ must occurs as a subquotient of $H^*_c(X_{\e^\chi, \s}(\e^\chi)$. In particular, it factors through $Z_{G(L), \s}(\e^\chi)/(Z_{G(L), \s}(\e^\chi) \cap I)$. See Corollary \ref{unique} and Corollary \ref{ddd}. 

(3) Let $G=PGL_n$ and $\t \in \tW$ is a superbasic element. Let $w \in \tW$. Then any simple representation of $Z_{G(L), \s}(\dot \t)$ occurs as a subquotient of $H^*_c(X_{w, \s}(\dot \t))$ is trivial. See Corollary \ref{superbasic}. 

\subsection{} We now review the content of this paper in more detail.

In section 1, we define good elements in extended affine Weyl groups and discuss some properties. In section 2, we discuss some reductive method and the relations between the dimension of affine Deligne-Lusztig varieties and the degrees of class polynomials for extended affine Hecke algebras. In section 3, we recall some combinatorial properties of extended affine Weyl group of type $A$ and $C$ established in \cite{He}. In section 4, we prove the main results and discuss some applications. 

\section{Good elements in extended affine Weyl group}

\subsection{} Let $G$ be a connected reductive algebraic group over an algebraically closed field $\kk$. Let $B$ be a Borel subgroup of $G$ and $B^-$ be an opposite Borel subgroup. Let $T=B \cap B^-$ be a maximal torus of $G$. Let $W_f=N_G(T)/T$ be the finite Weyl group of $G$. For $w \in W$, we choose a representative $\dot w \in N_G(T)$.

Let $G(L)$ be a loop group and $$\tW=X_*(T) \rtimes W_f=\{\e^\chi w; \chi \in X_*(T), w \in W_f\}$$ be the extended affine Weyl group of $G(L)$, where $\e$ is a symbol. The multiplication is given by the formula $(\e^\chi w) (\e^{\chi'} w')=\e^{\chi+w \chi'} w w'$. Let $I$ be the inverse image of $B^-$ under the projection map $G(o) \mapsto G$ sending $t$ to $0$ and let $I'$ be the inverse image of $U^-$, where $U^-$ is the unipotent radical of $B^-$. The we have the Bruhat-Tits decomposition $G(L)=\sqcup_{w \in \tW} I \dot{w} I$. It is also easy to see that if $\t \in \tW$ with $l(\t)=0$, then $\dot \t I \dot \t \i=I$. For $x=\e^\chi w \in \tW$, we choose a representative $\dot x=\e^\chi \dot w$ in $G(L)$. 

\subsection{} Let $R$ be the set of roots of $G$ and $R^+$ (resp. $R^-$) be the set of positive (resp. negative) roots of $G$. Let $(\a_i)_{i \in S}$ be the set of simple roots of $G$. For any $i \in S$, let $s_i$ be the corresponding simple reflection in $W_f$. Set $\tS=S \cup \{0\}$ and $s_0=t^{\th^\vee} s_\th$, where $\th$ is the largest positive root of $G$. 

Let $Q^\vee$ be the coroot lattice of $G$. Set $W_a=Q^\vee \rtimes W_f$. Then it is known that $W_a$ is a Coxeter group  with generators $s_i$ (for $i \in \tS$) and is a normal subgroup of $\tW$. Following \cite{IM65}, we define the length function on $\tW$ by $$l(\e^\chi w)=\sum_{\a \in R^+, w \i(\a) \in R^+} |<\chi, \a>|+\sum_{\a \in R^+, w \i(\a) \in R^-} |<\chi, \a>-1|.$$ It is known that for $w \in W_a$, $l(w)$ is the length of any reduced expression of $w$. For any coset of $W_a$ in $\tW$, there is a unique element of length $0$. Moreover, there is a natural group isomorphism between $\{\t \in \tW; l(\t)=0\}$ and $\tW/W_a \cong X_*(T)/Q^\vee$. 

Let $\t \in \tW$ with $l(\t)=0$, then for any $w, w' \in W_a$, we say that $\t w \le \t w'$ if $w \le w'$ for the Bruhat order on $W_a$. 

For any $J \subsetneqq \tS$, let $W_J$ be the subgroup of $\tW$ generated by $s_j$ (for $j \in J$), $w_0^J$ be the longest element in $W_J$ and ${}^J \tW$ be the set of minimal elements for the cosets $W_J \backslash \tW$. For $w \in {}^S \tW$, set $$I(w)=\max\{J \subset S; \forall j \in J, \exists j' \in J, \text{ such that } s_j w=w s_{j'}\}.$$ For any dominant coweight $\chi$, set $I(\chi)=\{i \in S; <\chi, \a_i>=0\}$. Then $I(\chi)=I(\e^\chi)$. 

For any subset $C$ of $\tW$, we set \begin{gather*} C_{\min}=\{w \in C; l(w) \le l(w') \text{ for any } w' \in C\}, \\ l(C)=l(w) \qquad \text{ for any } w \in C_{\min}. \end{gather*}

\subsection{}\label{comm} For any $\a \in R$, let $u_\a: \kk \to G$ with $t u_\a(k) t\i=u_\a(\a(t) k)$ for $k \in \kk$. It is easy to see that $u_\a$ extends in a natural way to a homomorphism $L \to G(L)$. 

%Then for $\a, \b \in R$ with $\a \neq \pm \b$, there exists $c_{\a, \b, i, j} \in \kk$ such that $$(u_\a(x), u_\b(y))=\prod_{i \a+j \b \in R, i, j>0} u_{i \a+j \b}(c_{\a, \b, i, j} x^i y^j)$$ for all $x, y \in \kk$. See \cite[Proposition 8.2.3]{Spr}. 

Let $\tilde R=\{\a+n \d; \a \in R, n \in \ZZ\}$ be the set of real affine roots and $\tilde R^+=\{\a+n \d; \a \in R^+, n>0\}\sqcup \{\a+n \d; \a \in R^-, n \ge 0\}$ be the set of positive real affine roots. The affine simple roots are $-\a_i$ for $i \in S$ and $\a_0=\th+\d$. Then any positive real affine root $\a$ can be written in a unique way as $\sum_{i \in S} -a_i \a_i+a_0 \a_0$ with $a_i \in \NN \cup \{0\}$ for $i \in \tS$. We set $ht(\a)=\sum_{i \in \tS} a_i$. For any real root $\a+n \d$, we define $x_{\a+n \d}: \kk \to G(L)$ by $x_{\a+n \d}(k)=u_\a(k \e^n)$ for $k \in \kk$. 

Let $I_1$ be the inverse image of $U^-$ under the projection $G(o) \to G$. Then $I_1$ is generated by $x_{\a}(\kk)$ for $\a \in \tilde R^+$ and $I_1 \cap T(L)$. For $n>1$, let $I_n$ be the subgroup of $I$ generated by $x_{\a}(\kk)$ with $ht(\a) \ge n$ and $I_1 \cap T(L)$. Then it is easy to see that $I_n$ is a normal subgroup of $I$ for $n \in \NN$ and for $\a, \b \in \tilde R^+$ and $a, b \in \kk$, \[\tag{*} x_\a(a) x_\b (b) \in x_\b(b) x_\a(a) I_{ht(\a)+ht(\b)}.\] 

\subsection{}\label{ne} Notice that each conjugacy class of $\tW$ lies in a coset of $W_a$. Let $\eta: \tW \to \tW/W_a$ be the natural projection. Then $\eta$ is constant on each conjugacy class of $\tW$. 

Let ${X_*(T)}_{\QQ}=X_*(T) \otimes_\ZZ \QQ$. Then the action of $W_f$ on $X_*(T)$ extends in a natural way to an action on ${X_*(T)}_\QQ$ and the quotient ${X_*(T)}_\QQ/W_f$ can be identified with $${X_*(T)}_\QQ^+=\{\chi \in {X_*(T)}_\QQ; <\chi, \a> \ge 0, \forall \a \in R^+\}.$$ 

For each element $w \in \tW$, there exists $n \in \NN$ such that $w^n=\e^\chi$ for some $\chi \in X_*(T)$. Let $v_{w}=\chi/n \in {X_*(T)}_\QQ$ and $[v_{w}]$ the corresponding element in ${X_*(T)}_\QQ/W_f$. It is easy to see that $v_{w}$ is independent of the choice of $n$. Moreover, if $w, w'$ are in the same conjugacy class of $\tW$, then $w^n$ is conjugated to $(w')^n$. Hence $(w')^n=\e^{\chi'}$ for some $\mu' \in W_f \cdot \chi$. Therefore $v_{w}$ and $v_{w'}$ are in the same $W_f$-orbit. We call the map $w \mapsto [v_{w}]$ the {\it Newton map}. Then the Newton map is constant on each conjugacy class of $\tW$. 

Define $f: \tW \to {X_*(T)}^+_\QQ \times \tW/W_a$ by $w \mapsto ([v_{w}], \eta(w))$. Then the map $f$ is constant on each conjugacy class of $\tW$. This map is the restriction to $\tW$ of the map $G(L) \to {X_*(T)}^+_\QQ \times \tW/W_a$ in \cite[4.13]{Ko97}. We denote the image of the map $f$ by $B(\tW)$. 

\begin{lem}
Let $w \in \tW$. Then the following conditions are equivalent:

(1) For any $n \in \NN$, $l(w^n)=n l(w)$. 

(2) Let $v$ be the unique element in $W_f \cdot v_{w} \cap {X_*(T)}_\QQ^+$. Then $$l(w)=<v, 2 \rho>,$$ where $\rho$ is the half sum of positive roots. 
\end{lem}

Proof. Assume that $w^m=\e^\chi$ for some $\chi \in X_*(T)$. 

If $l(w^n)=n l(w)$ for any $n \in \NN$, then $l(w)=\frac{1}{m} l(\e^\chi)$. Let $\chi' \in W_f \cdot \chi \cap X_*(T)^+$, then $l(\e^\chi)=l(\e^{\chi'})=<\chi', 2 \rho>$ and $v=\chi'/m$. Hence $l(w)=\frac{1}{m}<\chi', 2 \rho>=<v, 2 \rho>$. 

On the other hand, if $l(w)=<v, 2 \rho>$, then $$l(w^m)=l(\e^\chi)=<\chi', 2 \rho>=m l(w).$$ Let $n \in \NN$, then there exists $k \in \NN$ such that $n \le k m$. We have that $l(w^{m k})=l(\e^{k \chi})=<k \chi', 2 \rho>=m k l(w)$ and $$m k l(w)=l(w^{m k}) \le l(w^n)+l(w^{m k-n}) \le n l(w)+(m k-n) l(w)=m k l(w).$$ Therefore, both inequalities above are actually equalities. In particular, $l(w^n)=n l(w)$. \qed

\subsection{} We call an element $w \in \tW$ a {\it good element} if it satisfies the conditions in the previous lemma. The following result characterizes the good elements.

\begin{prop}\label{good}
Let $C$ be a fiber of $f: \tW \to B(\tW)$ and $w \in C$. Then $w$ is a good element if and only if $w \in C_{\min}$. 
\end{prop}

Proof. Notice that for any $x \in \tW$ and $n \in \NN$, $l(x) \ge \frac{1}{n} l(x^n)$. In particular, $l(x) \ge <v, 2 \rho>$, where $v$ is the unique element in $W_f \cdot v_{x} \cap {X_*(T)}_\QQ^+$. Hence if $w$ is a good element, then $w$ is a minimal length element in $C$.

Let $O$ be the $\s_F$-conjugacy class on $G(L)$ whose image under the map $G(L) \to {X_*(T)}^+_\QQ \times \tW/W_a$ equals $f(C)$. By \cite[Proposition 13.1.3 \& Corollary 13.2.4]{GHKR}, there exists a good element $x$ such that $I \dot x I \subset O$. Since $\dot x \in O$, $x \in C$. Therefore $x$ is a minimal length element and all the minimal length elements in $C$ are good elements.  \qed

\begin{cor}
Let $w$ be a Coxeter element in $W_a$. Then $w$ is a good element. 
\end{cor}

\begin{rmk}
This result was first proved by Speyer in \cite{Spe}. Here we give a different proof. 
\end{rmk}

Proof. By \cite{Ho82}, $w$ has infinite order. Hence $[w] \neq 0$. Let $C$ be a fiber of the map $f: \tW \to B(\tW)$ that contains $w$. If $w' \in C$ and $l(w')<l(w)$, then $w'$ lies in some $W_J$ with $J \neq \tS$. In particular, $w'$ lies in some finite Weyl group and $[w']=0$. That is a contradiction. So $w$ is a minimal length element in $C$. By the previous Proposition, $w$ is a good element. \qed

\begin{lem}\label{dg}
Let $w, w'$ be good elements. If $f(w) \neq f(w')$, then $$G(L) \cdot_\s I \dot w I \cap G(L) \cdot_\s I \dot w' I=\emptyset.$$ 
\end{lem}

Proof. It is easy to see that for $x,\ty \in \tW$, $I \dot x I \dot \ty I \subset \cup_{\tz \in x \ty W_a} I \tz I$. Hence \begin{gather*} G(L) \cdot_\s I \dot w I \subset \cup_{x \in \tW} I \dot x I \dot w I \dot x \i I \subset \cup_{\tz \in x w x \i W_a=w W_a} I \dot \tz I, \\ G(L) \cdot_\s I \dot w' I \subset \cup_{x \in \tW} I \dot x I \dot w' I \dot x \i I \subset \cup_{\tz \in x w' x \i W_a=w W_a} I \dot \tz I. \end{gather*} Therefore if $w W_a \neq w' W_a$, then $G(L) \cdot_\s I \dot w I \cap G(L) \cdot_\s I \dot w' I=\emptyset$. 

Now assume that $\eta(w)=\eta(w')$. By our assumption, $[v_{w}] \neq [v_{w'}]$. If $G(L) \cdot_\s I \dot w I \cap G(L) \cdot_\s I \dot w' I \neq \emptyset$, then there exist $g \in G(L)$ such that $I \dot w I \cap g I \dot w' I \s(g) \i \neq \emptyset$. Let $z \in I \dot w I \cap g I \dot w' I \s(g) \i$. Since $w$ and $w'$ are good elements, then for any $n \in \NN$, \begin{align*} z \s(z) \cdots \s^{n-1}(z) & \in (I \dot w I) (I \dot w I) \cdots (I \dot w I)=I \dot w^n I, \\ z \s(z) \cdots \s^{n-1}(z) & \in (g I \dot w' I \s(g) \i) (\s(g) I \dot w' I \s^2(g) \i) \cdots (\s^{n-1}(g) I \dot w' I \s^n(g) \i) \\ &=g I (\dot w')^n I \s^n(g). \end{align*} In particular, for any $n \in \NN$, $I \dot w^n I \cap g I (\dot w')^n I \s^n(g) \i \neq \emptyset$. 

There exists $m \in \NN$ such that $w^m=\e^\mu$ and $(w')^m=\e^{\mu'}$ for some $\mu, \mu' \in X_*(T)$. Since $[v_{w}] \neq [v_{w'}]$, $\mu \notin W_f \cdot \mu'$. Assume that $g \in I x I$ for some $x \in \tW$. Then $I \dot \e^{k \mu} I \cap I \dot x I \dot \e^{k \mu'} I \dot x \i I \neq \emptyset$ for all $k \in \NN$. 

Notice that $I \dot x I \dot \e^{k \mu'} I \dot x \i I \subset \cup_{\ty, \ty' \le x} I \dot \ty \dot \e^{k \mu'} (\dot \ty') \i I$. Thus $\e^{k \mu}=\ty \e^{k \mu'} (\ty') \i$ for some $\ty, \ty' \le x$. Assume that $\ty=y \e^\chi$ and $\ty'=y' \e^{\chi'}$ with $\chi, \chi' \in X_*(T)$ and $y, y' \in W_f$. Then $\ty \e^{k \mu'} (\ty') \i=\e^{y (k \mu'+\chi-\chi')} y (y')\i$. Hence $y=y'$ and $k \mu'+\chi-\chi'=k y \i \mu$. By definition, $l(\e^\chi) \le l(\ty)+l(y) \le l(x)+l(w_0^S)$. Similarly, $l(\e^{\chi'}) \le l(x)+l(w_0)$. Since $\mu \notin W_f \cdot \mu'$, then $l(\e^{\mu'-y \i \mu}) \ge 1$. Now $$k \le l(\e^{k(\mu'-y \i \mu)})=l(\e^{\chi'-\chi}) \le l(\e^{\chi'})+l(\e^\chi) \le 2 l(x)+2 l(w_0).$$ That is a contradiction. \qed 

\begin{cor}\label{single}
Let $w, w' \in \tW$ are good elements with $f(w)=f(w')$. Then $G(L) \cdot_{\s_F} I \dot w I=G(L) \cdot_{\s_F} I \dot w' I$ is a single $\s_F$-conjugacy class. 
\end{cor}

Proof. As in the proof of Proposition \ref{good}, for any $\s_F$-conjugacy class $O$ of $G(L)$, there exists a good element $x_O \in \tW$ such that $\dot x_O \in O$. If $O \subset G(L) \cdot_{\s_F} I \dot w I$, then we must have that $\dot x_O \in G(L) \cdot_{\s_F} I \dot w I$. By the previous Lemma, $f(x_O)=f(w)$. By \cite[4.13]{Ko97}, $O$ is uniquely determined by $f(x_O)$. Therefore $G(L) \cdot_{\s_F} I \dot w I$ is the single $\s_F$-conjugacy class that contains $\dot x_O$. In particular,  $G(L) \cdot_{\s_F} I \dot w I=G(L) \cdot_{\s_F} I \dot w' I$. \qed

\subsection{}\label{reK} Now we reformulate Kottwitz's classification of $\s_F$-conjugacy classes as follows.

Let $\tW_{good}$ be the set of good elements in $\tW$. For $w, w' \in \tW_{good}$, we write $w \sim w'$ if $f(w)=f(w')$. Then we have that \[\tag{*} G(L)=\sqcup_{[w] \in \tW_{good}/\sim} G(L) \cdot_{\s_F} I \dot w I.\]

However, if $\s=\s_a$, then in general, $G(L) \cdot_\s I \dot w I$ contains infinitely many $\s$-conjugacy classes. In that case, we don't know if $G(L) \cdot_\s I \dot w I=G(L) \cdot_\s I \dot w' I$ for $w \sim w'$. However, for some groups it still holds when $\s=\s_a$ and we have similar decomposition. We will discuss it in more details in section 4. 

\begin{lem}\label{0}
Let $b \in I$. Then $\dim(X_{1, \s_F}(b))=0$ and $\dim(X_{1, \s_a}(b))<\infty$ if $a$ is not a root of unity. 
\end{lem} 

Proof. By \cite[Prop 6.3.1]{GHKR}, any element in $I$ is $\s_F$ conjugate to $1$. So we may take $b=1$ if $\s=\s_F$. 

For $x \in \tW$, define $Y_x=I \dot x I/I \cap X_{1, \s}(b)$.  Assume that $\{\a \in \tilde R^+; x \i \a<0\}=\{\a_{i_1}, \cdots, \a_{i_k}\}$. We arrange the roots such that $ht(\a_{i_1}) \le ht(\a_{i_2}) \le \cdots \le ht(\a_{i_k})$. Let $g \in G(L)$ with $g I \in Y_x$, then $g=x_{\a_{i_1}}(a_1) \cdots x_{\a_{i_k}}(a_k) \dot x I$ for some $a_1, \cdots, a_k \in \kk$. Then $g I \cap b \s(g) I \neq \emptyset$. If $\s=\s_F$ and $b=1$, then we must have that $x_{\a_{i_1}}(a_1) \cdots x_{\a_{i_k}}(a_k)=x_{\a_{i_1}}(F(a_1)) \cdots x_{\a_{i_k}}(F(a_k))$. Therefore $a_1=F(a_1), \cdots, a_k=F(a_k)$. So there are only finite many elements in $Y_x$ for each $x \in \tW$ and $\dim(X_{1, \s_F}(b))=0$. 

Now we consider the case where $\s=\s_a$ with $a$ not a root of unity.  Notice that $\dim Y_x$ equals the dimension of the variety consists of $(a_1, \cdots, a_k)$ satisfying $$x_{\a_{i_1}}(a_1) \cdots x_{\a_{i_k}}(a_k) \dot x I \cap b \s(x_{\a_{i_1}}(a_1) \cdots x_{\a_{i_k}}(a_k)) \dot x I \neq \emptyset.$$  

We may assume that $b \in t I_1$. Let $n \in \NN$ with $t \s(x_\a(1)) t \i \neq x_\a(1)$ for all $\a$ with $ht(\a) \ge n$. We show that $$\dim(Y_x) \le \sharp\{\a \in \tilde R^+; ht(\a)<n\}$$ for any $x \in \tW$. 

We may assume that $ht(\a_{i_{j-1}})<n$ and $ht(\a_{i_j})=n$. It is enough to show that for any $u \in I_n ({}^{\dot x} I \cap I)$, there is a unique $(a_j, \cdots, a_k)$ such that $x_{\a_{i_j}}(a_j) \cdots x_{\a_{i_k}}(a_k) \in u t \s(x_{\a_{i_j}}(a_j) \cdots x_{\a_{i_k}}(a_k)) ({}^{\dot x} I \cap I)$. 

We prove this statement by descending induction on $n$. 

Assume that $ht(\a_{i_j})=\cdots=ht(\a_{i_l})=n<ht(\a_{i_{l+1}})$. Then \begin{gather*} x_{\a_{i_j}}(a_j) \cdots x_{\a_{i_k}}(a_k) \in x_{\a_{i_j}}(a_j) \cdots x_{\a_{i_l}}(a_l) I_{n+1}, \\ u t \s(x_{\a_{i_j}}(a_j) \cdots x_{\a_{i_k}}(a_k)) \in u' t \s(x_{\a_{i_j}}(a_j) \cdots x_{\a_{i_l}}(a_l)) I_{n+1} ({}^{\dot x} I \cap I),\end{gather*} where $u'$ is the unique element in $x_{\a_{i_j}}(\kk) \cdots x_{\a_{i_l}}(\kk) \cap u I_{n+1} ({}^{\dot x} I \cap I)$. Since $[{}^{\dot x} I \cap I, I_n] \subset I_{n+1}$, $$x_{\a_{i_j}}(a_j) \cdots x_{\a_{i_l}}(a_l) \in u' t \s(x_{\a_{i_j}}(a_j) \cdots x_{\a_{i_l}}(a_l)) t\i I_{n+1} ({}^{\dot x} I \cap I).$$ Assume that $u'=x_{\a_{i_j}}(b_j) \cdots x_{\a_{i_l}}(b_l)$. For $p \ge j$, $t \s(x_{\a_p}(k)) t \i=x_{\a_p}(c_p k)$ for some $c_p \neq 1$. Thus we must have that $a_j=b_j+c_j a_j, \cdots, a_l=b_l+c_l a_l$. In particular, $a_j, \cdots, a_l$ are uniquely determined. Now set $u_1=(x_{\a_{i_j}}(a_j) \cdots x_{\a_{i_l}}(a_l)) \i u' t \s(x_{\a_{i_j}}(a_j) \cdots x_{\a_{i_l}}(a_l)) t \i \in I_{n+1} ({}^{\dot x} I \cap I)$, we have that $$x_{\a_{i_{l+1}}}(a_{l+1}) \cdots x_{\a_{i_k}}(a_k) \in u_1 t \s(x_{\a_{i_{l+1}}}(a_{l+1}) \cdots x_{\a_{i_k}}(a_k)) ({}^{\dot x} I \cap I).$$ By induction hypothesis, $a_{l+1}, \cdots, a_k$ are also uniquely determined. \qed

\begin{lem}\label{finite}
Let $\mu \in X_*(T)$ and $M$ be the Levi subgroup of $G$ generated by $T$ and $u_\a(\kk)$ for $\a \in R$ with $<\mu, \a>=0$. Then the map $$I \times_{I \cap M(L)} (I \cap M(L)) \e^\mu \to I \e^\mu I$$ defined by $(i, i') \mapsto i i' \s(i) \i$, is bijective. 
\end{lem}

\begin{rmk}
If $\s=\s_F$, then the lemma is a special case of \cite[Theorem 2.1.2]{GHKR}. The case where $\s=\s_a$ is essentially the same as in loc.cit. Here we give a proof to convince the readers that no problem occurs for $\s=\s_a$. 
\end{rmk}

Proof. Notice that $I \times_{I \cap M(L)} (I \cap M(L)) \e^\mu \cong I_1 \times_{I_1 \cap M(L)} (I \cap M(L)) \e^\mu$ and $I \e^\mu I=I_1 (I \cap M(L) \e^\mu I_1$. It suffices to prove that for any $n$, the map $I_n \times_{I_n \cap I_{n+1} (I_n \cap M(L))} I_{n+1} (I \cap M(L)) \e^\mu I_{n+1} \to I_n (I \cap M(L)) \e^\mu I_n$ defined by $(i, i') \mapsto i i' \s(i) \i$ is bijective. 

Let $P$ be the parabolic subgroup of $G$ generated by $T$ and $u_\a(\kk)$ for $\a \in R$ with $<\mu, \a> \ge 0$ and Let $P^-$ be the opposite parabolic subgroup of $G$ generated by $T$ and $u_\a(\kk)$ for $\a \in R$ with $<\mu, \a> \le 0$. Set $I_n'=I_n \cap U_P(L)$ and $I_n''=I_n \cap U_{P^-}(L)$. Since $I$ normalizes $I_n$ and $M(L)$ normalizes $U_{P}(L)$ and $U_{P^-}(L)$, $I \cap M(L)$ normalizes $I'_n$ and $I''_n$ for all $n$. Moreover, we have that $I_n=I'_n I''_n (I_n \cap M(L))$. Now \begin{align*} I_n (I \cap M(L)) \e^\mu I_n & \subset I_n \cdot_\s I_n (I \cap M(L)) \e^\mu=I_n \cdot_\s I'_n I''_n (I \cap M(L)) \e^\mu \\ &=I_n \cdot_\s I''_n (I \cap M(L)) \e^\mu I'_n=I_n \cdot_\s (I \cap M(L)) I''_n \e^\mu I'_n. \end{align*} By definition, $\e^\mu I'_n \e^{-\mu} \subset I_{n+1}$ and $\e^{-\mu} I''_n \e^\mu \subset I_{n+1}$. So $I''_n \e^\mu I'_n \subset I''_n I_{n+1} \e^\mu=I_{n+1} I''_n \e^\mu \subset I_{n+1} \e^\mu I_{n+1}$. This proves the surjectivity.  

On the other hand, for any $i_1 \in I'_n$ and $i_2 \in I''_n$, \begin{align*} & i_1 i_2 (I_n \cap M(L)) \cdot_\s I_{n+1} (I \cap M(L)) \e^\mu I_{n+1}=i_1 i_2 I_{n+1} (I \cap M(L)) \e^\mu I_{n+1} \s(i_1 i_2) \i  \\ & \subset i_1 I_{n+1} (I \cap M(L)) I''_n \e^\mu I'_n I_{n+1} \s(i_2) \i=i_1 I_{n+1} (I \cap M(L)) \e^\mu I_{n+1} \s(i_2) \i.\end{align*} So $i_1 i_2 (I_n \cap M(L)) \cdot_\s I_{n+1} (I \cap M(L)) \e^\mu I_{n+1} \cap I_{n+1} (I \cap M(L)) \e^\mu I_{n+1} \neq \emptyset$ if and only if $i_1, i_2 \in I_{n+1}$. This proves the injectivity. \qed

\begin{prop}\label{00}
Let $w \in \tW$ be a good element. Let $n \in \NN$ and $\mu \in X_*(T)$ with $w^n=\e^\mu$. Let $M$ be the Levi subgroup of $G$ generated by $T$ and $u_\a(\kk)$ for $\a \in R$ with $<\mu, \a>=0$. Then $X_{w, \s}(\dot w) \subset M(L) I/I \cong M(L)/(M(L) \cap I)$. Moreover, for any $b \in G(L) \cdot_\s I \dot w I$, $\dim(X_{w, \s}(b))=0$ if $\s=\s_F$ and $\dim(X_{w, \s}(b))<\infty$ if $\s=\s_a$ with $a$ not root of unity. 
\end{prop}

Proof. Since $w$ is a good element, $(I \dot w I) \s(I \dot w I) \cdots \s^{n-1}(I \dot w I) \subset I \e^\mu I$. By Lemma \ref{finite}, for any $b \in G(L) \cdot_\s I \dot w I$, there exists $h \in G(L)$ such that $h \i b \s(b) \cdots \s^{n-1}(b) \s^n(h) \in (I \cap M(L)) \e^\mu$. If $b=\dot w$, then we may just take $h=1$. Now set $b'=h \i b \s(b) \cdots \s^{n-1}(b) \s^n(h)$. If $g I \in X_{w, \s}(b)$, then again by Lemma \ref{finite}, $g I=g' I$ for some $g'$ with $(g') \i b \s(b) \cdots \s^{n-1}(b) \s^n(g') \in (I \cap M(L)) \e^\mu$. Thus $$X_{w, \s}(b) \subset \{ h g I/I; g \i b' \s^n(g) \in (I \cap M(L)) \e^\mu\}.$$ 

Let $x \in W_f$ such that $x \mu$ is dominant. Set $J=I(x \mu)$ and $M'={}^{\dot x} M$. Then $M'$ is a Levi factor of $P=P_J$ and ${}^{\dot x} (I \cap M(L)) \subset \cap M'(o)$. It is easy to see that $$X_{w, \s}(b) \subset \{ h \dot x g \dot x \i I/I; g \i M'(o) \e^{x \mu} \s^n(g) \cap M'(o) \e^{x \mu} \neq \emptyset\}.$$

Let $g \in G(L)$ with $g \i M'(o) \e^{x \mu} \s^n(g) \cap M'(o) \e^{x \mu} \neq \emptyset$. Then we must have that $g=m k$ for $m \in M'(L)$ and $k \in G(o)$ and $m \i M'(o) \e^{x \mu} \s^n(m) \in M'(o) \e^{x \mu}$. The case where $\s=\s_F$ is proved in \cite[Theorem 1.1 (2)]{Ko2}. The case where $\s=\s_a$ can be proved in the same way. 

Assume that $k \in I \dot u I$ for $u=u_1 v_1$ with $u_1 \in W^J$ and $v_1 \in W_J$. Then $I \dot u I=I \dot u_1 I \dot v_1 I \subset I \dot u_1 I M'$. Thus $I \dot u_1 I M'(o) \e^{x \mu} \cap M'(o) \e^{x \mu} \s^n(I \dot u_1 I) \neq \emptyset$. Notice that \begin{gather*} M'(o) \e^{x \mu} \s^n(I \dot u_1 I) \subset \sqcup_{v \in W_J} I \dot v I \e^{x \mu} I \dot u_1 I=\sqcup_{v \in W_J} I \e^{x \mu} \dot v I \dot u_1 I \subset \sqcup_{y \in W_f} I \e^{x \mu} \dot y I, \\ I \dot u_1 I M'(o) \e^{x \mu} \subset \sqcup_{v \in W_J} I \dot u_1 I \dot v I \e^{x \mu}=\sqcup_{v \in W_J} I \dot u_1 I \e^{x \mu} I \dot v I=\sqcup_{v \in W_J} I \dot u_1 \e^{x \mu} \dot v I. \end{gather*} Hence there exists $v \in W_J$ and $y \in W$ such that $\e^{x \mu} y=u_1 \e^{x \mu} v=\e^{u_1 x \mu} u_1 v$. In particular, we have that $x \mu=u_1 x \mu$. By \cite[Lemma 3.5]{HT}, $u_1 \in W_J$. So $u_1=1$ and $k \in I M'(o)$. Since $k \i M'(o) \e^{x \mu} \s^n(k) \cap M'(o) \e^{x \mu} \neq \emptyset$, by the proof of Lemma \ref{finite}, $k \in M'(o)$. Therefore $X_{w, \s}(b) \subset \{h \dot x g \dot x \i I/I; g \in M'(L)\}=\{h g I/I; g \in M(L)\}$. The ``moreover'' part follows from Lemma \ref{0}.

\section{Reductive method}

\subsection{} Let $G(L)'=\sqcup_{w \in W_a} I \dot w I$ be the identity component of $G(L)$. Then \begin{align*} G(L) &=\sqcup_{w \in \tW} I \dot w I=\sqcup_{\t \in \tW, l(\t)=0} \sqcup_{w \in W_a} I \dot w \dot \t I \\&=\sqcup_{\t \in \tW, l(\t)=0} \sqcup_{w \in W_a} I \dot w I \dot \t=\sqcup_{\t \in \tW, l(\t)=0} G(L)' \dot \t. \end{align*}

\subsection{} We follow the notations in \cite[2.1]{He}. 

For $w, w' \in \tW$ and $i \in \tS$, we write $w \xrightarrow{s_i} w'$ if $w'=s_i w s_i$ and $l(w') \le l(w)$. We write $w \to w'$ if there is a sequence $w=w_0, w_1, \cdots, w_n=w'$ of elements in $\tW$ such that for all $k$, $w_{k-1} \xrightarrow{s_i} w_k$ for some $i \in \tS$. We write $w \tilde \to w'$ if there exists $\t \in \tW$ with $l(\t)=0$ such that $w \to \t w' \t \i$, or equivalently, there is a sequence $w=w_0, w_1, \cdots, w_n=w'$ of elements in $\tW$ such that for all $k$, $w_k=\t w_{k-1} \t \i$ for some $\t \in \tW$ with $l(\t)=0$ or $w_{k-1} \xrightarrow{s_i} w_k$ for some $i \in \tS$. We write $w \tilde \approx w'$ if $w \tilde \to w'$ and $w' \tilde \to w$. We write $w \approx w'$ if $w \to w'$ and $w' \to w$.

\begin{lem}\label{red1}
Let $w, w' \in \tW$. 

(1) If $w \to w'$, then $$G(L)' \cdot_\s I \dot w I \subset G(L)' \cdot_\s I \dot w' I \cup \cup_{x \in w W_a, l(x)<l(w)} G(L)' \cdot_\s I \dot x I.$$ 

(2) If $w \approx w'$, then $$G(L)' \cdot_\s I \dot w I=G(L)' \cdot_\s I \dot w' I.$$ 
\end{lem}

Proof. By definition, there exists a finite sequence $w=w_0 \xrightarrow {i_1} w_1 \xrightarrow {i_2} \cdots \xrightarrow {i_m} w_n=w'$, where $i_j \in \tS$ for all $j$. We prove the lemma by induction on $m$. 

The statements are true for $m=0$. Now assume that the statements hold for $m-1$. By \cite[Lemma 1.6.4]{DL}, we have that $w=w_1$ or $s_{i_1} w<w$ or $w s_{i_1}<w$. If $w=w_1$, then the statements follow from induction hypothesis. Now we prove the case where $s_{i_1} w<w$. The case $w s_{i_1}<w$ can be proved in the same way. 

Since $s_{i_1} w<w$, then $G(L)' \cdot_\s I \dot w I=G(L)' \cdot_\s I \dot s_{i_1} I \dot s_{i_1} \dot w I=G(L)' \cdot_\s I \dot s_{i_1} \dot w I \dot s_{i_1} I$. Moreover, \[I \dot s_{i_1} \dot w I \dot s_{i_1} I=\begin{cases} I \dot w_1 I, & \text{ if }  l(w_1)=l(s_{i_1} w)+1=l(w); \\ I \dot s_{i_1} \dot w I \sqcup I \dot w_1 I, & \text{ if } l(w_1)=l(s_i w)-1=l(w)-2. \end{cases}\] In either case, $$G(L)' \cdot_\s I \dot w I \subset G(L)' \cdot_\s I \dot w_1 I \cup \cup_{x \in w W_a, l(x)<l(w)} G(L)' \cdot_\s I \dot x I.$$ Notice that $l(w_1) \le l(w)$. By induction hypothesis, $G(L)' \cdot_\s I \dot w_1 I \subset G(L)' \cdot_\s I \dot w' I \cup \cup_{x \in w W_a, l(x)<l(w_1)} G(L)' \cdot_\s I \dot x I$. Hence $G(L)' \cdot_\s I \dot w I \subset G(L)' \cdot_\s I \dot w' I \cup \cup_{x \in w W_a, l(x)<l(w)} G(L)' \cdot_\s I \dot x I$. 

If moreover, $w \approx w'$, then $l(w_1)=l(w)$ and $w_1 \approx w'$. By induction hypothesis, $G(L)' \cdot_\s I \dot w I=G(L)' \cdot_\s I \dot w_1 I=G(L)' \cdot_\s I \dot w' I$. \qed

\begin{lem}\label{red2}
Let $w, w' \in \tW$. 

(1) If $w \tilde \to w'$, then $$G(L) \cdot_\s I \dot w I \subset G(L) \cdot_\s I \dot w' I \cup \cup_{x \in w W_a, l(x)<l(w)} G(L) \cdot_\s I \dot x I.$$ 

(2) If $w \tilde \approx w'$, then $$G(L) \cdot_\s I \dot w I=G(L) \cdot_\s I \dot w' I.$$ 
\end{lem}

Proof. For any $x \in \tW$ and $\t \in \tW$ with $l(\t)=0$, we have that $$G(L) \cdot_\s I \dot x I=G(L) \cdot_\s \dot \t I \dot x I \s(\dot \t) \i=G(L) \cdot_\s I \dot \t \dot x \dot \t \i I.$$ Now if $w \tilde \to w'$, then there exists $\t \in \tW$ with $l(\t)=0$ such that $w \to \t w' \t \i$. By the previous lemma,  \begin{align*} G(L) \cdot_\s I \dot w I &=G(L) \cdot_\s (G(L)' \cdot_\s I \dot w I) \\ &\subset G(L) \cdot_\s \bigl(G(L)' \cdot_\s I \dot \t \dot w' \dot \t \i I \cup \cup_{x \in w W_a, l(x)<l(w)} G(L)' \cdot_\s I \dot x I \bigr) \\ &=G(L) \cdot_\s I \dot \t \dot w' \t \dot \i I \cup \cup_{l(x)<l(w)} G(L) \cdot_\s I \dot x I \\ &=G(L) \cdot_\s I \dot w' I \cup \cup_{x \in w W_a, l(x)<l(w)} G(L) \cdot_\s I \dot x I. \end{align*} 

If $w \tilde \approx w'$, then there exists $\t \in \tW$ with $l(\t)=0$ such that $w \approx \t w' \t \i$. By the previous lemma,  \begin{align*} G(L) \cdot_\s I \dot w I &=G(L) \cdot_\s (G(L)' \cdot_\s I \dot w I) =G(L) \cdot_\s (G(L)' \cdot_\s I \dot \t \dot w' \dot \t \i I) \\ &=G(L) \cdot_\s I \dot \t \dot w' \dot \t \i I=G(L) \cdot_\s I \dot w' I. \end{align*} The Lemma is proved. \qed

\begin{lem}\label{red3} Let $w \in {}^S \tW$ and $w'=x w$ for some $x \in W_{I(w)}$, then $G(L)' \cdot_\s I \dot w' I \subset G(L)' \cdot_\s I \dot w I$. 
\end{lem}

Proof. Set $J=I(w)$. Notice that $I=(I \cap L_J) I_J$, where $I_J$ is the inverse image of $U_{P_J^-}$ under the map $G(o) \to G$. It is easy to see that $I_J$ is a normal subgroup of $I$ and $L_J$ normalizes $I_J$. Thus $I \dot w I=(I \cap L_J) I_J \dot w (I \cap L_J) I_J=(I \cap L_J) I_J (I \cap L_J) \dot w I_J=(I \cap L_J) I_J \dot w I_J$. We have that $$\cup_{x \in W_J} I \dot x \dot w I=\cup_{x \in W_J} B^- \dot x B^- I \dot w I=P_J^- I \dot w I=L_J I \dot w I=L_J I_J \dot w I_J.$$ Define $\s': L_J \to L_J$ by $\s'(l)=\dot w \s(l) \dot w \i$. Then for $l, l' \in L_J$, $$l (l' I_J \dot w I_J) \s(l) \i=l l' (\s' \circ \s)(l) \i I_J \dot w I_J.$$ Notice that if $\s$ is an automorphism on $L_J$, then $\s'$ is also an automorphism on $L_J$. If $\s$ is a Frobenius morphism on $L_J$, then $\s'$ is also a Frobenius morphism on $L_J$. Moreover, $\s'(B^- \cap L_J)=B^- \cap L_J$. By \cite[Lemma 7.3]{St68} (if $\s'$ is an automorphism) and Lang's theorem \cite{La56} (if $\s'$ is a Frobenius morphism), $L_J \cdot_{\s'} (B^- \cap L_J)=L_J$. Hence $L_J I_J \dot w I_J=(L_J) \cdot_\s (B^- \cap L_J) I_J \dot w I_J \subset (L_J) \cdot_\s I \dot w I$. \qed

\

Now we discuss some reductive method for $X_{w, \s}(b)$. The following result can be proved along the line of the proof of \cite[Theorem 1.6]{DL}. 

\begin{lem}\label{indd}
Let $x\in \widetilde{W}$, and let $s\in\tS$ be a simple affine reflection. Then 

(1) If $l(sxs)=l(x)$, then $X_{x, \s}(b)$, $X_{sxs, \s}(b)$ are universally homeomorphic. In this case, $H^*_c(X_{x, \s}(b)) \cong H^*_c(X_{s x s, \s}(b))$ as $Z_{G(L), \s}(b)$-modules. 

(2) If $l(sxs) =l(x)-2$, then $X_{x, \s}(b)$ can be written as a disjoint union $X_{x, \s}(b) = X_1 \cup X_2$ where $X_1$ is closed and $X_2$ is open, and such that $X_1$ admits a morphism to $X_{sxs, \s}(b)$, all of whose fibers are isomorphic to $\mathbb A^1$, and such that $X_2$ admits a morphism to $X_{sx, \s}(b)$, all of whose fibers are isomorphic to $\mathbb A^1\setminus\{ 0\}$. In this case, any simple $Z_{G(L), \s}(b)$-module occurs as a subquotient of $H^*_c(X_{x, \s}(b))$ must occurs as a subquotient of $H^*_c(X_{s x s, \s}(b))$ or $H^*_c(X_{s x, \s}(b))$. 
\end{lem}

We also have the following result. The case where $\s=\s_F$ is proved in \cite{GH}. The case where $\s=\s_a$ can be proved in the same way. 

\begin{lem}\label{tt}
Let $x, \t \in \tW$ with $l(\t)=0$. Then for any $b \in G(L)$, $X_{x, \s}(b)$ is isomorphic to $X_{\t x \t \i, \s}(b)$. 
\end{lem}

\begin{lem}\label{lang}
Let $w \in {}^S \tW$ and $x \in W_{I(w)}$. Then $\dim(X_{x w, \s_F}(b))=\dim(X_{w, \s_F}(b))+l(x)$ and $\dim(X_{x w, \s_a}(b)) \le \dim(X_{w, \s_a}(b))+l(w_0^{I(w)})$.
\end{lem}

\begin{rmk}
By convention, we set the dimension of the empty set to be $-\infty$. 
\end{rmk}

Proof. Set $J=I(w)$ and $P=L_J I$. Define $$X=\{g P; g \i b \s(g) \in P w P\}.$$ Notice that for any $u \in W_J$, $I \dot u \dot w I=I \dot u I \dot w I\subset I \dot w P \subset P \dot w P$. Thus the map $g I \mapsto g P$ sends $X_{u w, \s}(b)$ to $X$. 

Let $g P \in X$. Then $g \i b \s(g) \in P \dot w P=I_J L_J \dot w I_J$, where $I_J$ is the inverse image of $U_{P_J^-}$ under the map $G(o) \to G$. We assume that $g \i b \s(g) \in I_J l \dot w I_J$ for $l \in L_J$. Now for $p \in L_J$, $p \i g \i b \s(g) \s(p) \in p \i I_J l \dot w I_J \s(p) \i=I_J p \i l \dot w \s(p) \i I_J$. Notice that $I \dot u \dot w I=I_J (B^- \cap L_J) \dot u (B^- \cap L_J) \dot w I_J$. Thus for $p \in L_J$, $g p I \in X_{u w, \s}(b)$ if and only if $p \i l \dot w \s(p) \dot w \i \in (I \cap L_J) \dot u (I \cap L_J)$. Define $\s': L_J \to L_J$ by $\s'(l)=\dot w \s(l) \dot w \i$. Set $Y_g=\{p (I \cap L_J) \in L_J/(I \cap L_J); p \i l \s'(p) \in (B^- \cap L_J) \dot u (B^- \cap L_J)\}$. Then $\{p I \in P/I; g p I \in X_{u w, \s} (b)\} \cong Y_g$. Notice that $L_J/(I \cap L_J)=L_J/(L_J \cap B^-)$. So $Y_g$ is a subvariety of $L_J/(L_J \cap B^-)$. In particular, each fiber of the map $X_{u w, \s}(b) \to X$ is at most of dimension $l(w_0^J)$. 

If $\s=\s_F$, $Y_g$ is a Deligne-Lusztig variety in $L_J/(L_J \cap B^-)$ and it is known that it is of dimension $l(u)$. So $\dim(X_{u w, \s_F}(b))=\dim(X)+l(u)$ for any $u \in W_J$. In particular, $\dim(X_{w, \s_F}(b))=\dim(X)$ and $\dim(X_{x w, \s_F}(b))=\dim(X_{w, \s_F}(b))+l(x)$. 

If $\s=\s_a$, then by \cite[Lemma 7.3]{St68}, $Y_g$ is nonempty for $u=1$. Hence $\dim(X_{w, \s_a}(b)) \ge \dim(X)$. Therefore $\dim(X_{x w, \s_a}(b)) \le \dim(X)+l(w_0^J) \le \dim(X_{w, \s_a}(b))+l(w_0^J)$. \qed

\

Now let us make a short digression and discuss the class polynomials of extended affine Hecke algebra. We will then discuss the relation between the degree of these polynomials and the dimension of $X_{x, \s}(b)$ in the end of this section. 

\subsection{}\label{Hecke} Let $H$ be the Hecke algebra associated to an extended affine Weyl group $\tW$, i.e., $H$ is the associated $\ZZ[v, v \i]$-algebra with basis $T_w$ for $w \in \tW$ and multiplication is given by 
\begin{gather*} T_x T_y=T_{x y}, \quad \text{ if } l(x)+l(y)=l(x y); \\ (T_s-v)(T_s+v \i)=0, \quad \text{ for } s \in \tS. \end{gather*}

Then $T_s \i=T_s-(v-v \i)$ and $T_w$ is invertible in $H$ for all $w \in \tW$. 

If $w, w' \in \tW$ with $w'=x w x \i$, $l(w')=l(w)$ and $l(x w)=l(x)+l(w)$, then $T_x T_w=T_{x w}=T_{w' x}=T_{w'} T_x$ and $T_w-T_{w'}=T_x \i T_{x w}-T_{x w} T_x \i=[T_x \i, T_{x w}] \in [H, H]$. Therefore, we have that 

(1) If $w \tilde \sim w'$, then $T_w \equiv T_{w'} \mod [H, H]$. 

It is also easy to check that 

(2) If $l(s w s)<l(w)$ for some $s \in \tS$, then $$T_w \equiv T_s^2 T_{s w s}=(v-v \i) T_{w s}+T_{s w s} \mod [H, H].$$ 

\subsection{} In the rest of this section, we assume that $G$ is a simple algebraic group of type $A_{n-1}$, $B_n$, $C_n$ or $D_n$ and $\tW$ be the extended affine Weyl group for $G$. We embed the root lattice and coweight lattice in $\oplus_{i=1}^n \RR e_i$ in the natural way (see, e.g., \cite[Plate I-IV]{Bo}). An element $\e^\chi w \in \tW$ with $w \in W$ and $\chi \in \oplus_{i=1}^n \ZZ e_i$ is call an {\it integral element}. We denote by $\tW_{int}$ the subset of all integral elements in $\tW$. Set \begin{gather*} G(L)^!=\begin{cases} G(L) \rtimes <\iota>, &\text{ if } G \text{ is of type } D\\ G(L), & \text{ otherwise} \end{cases}; \\ \tW^!=\begin{cases} \tW \rtimes <\iota>, &\text{ if } G \text{ is of type } D\\ \tW, & \text{ otherwise} \end{cases}.\end{gather*} Here $\iota$ is the outer diagonal automorphism on $G$ whose induces action on $\oplus_{i=1}^n \RR e_i$ sends $e_n$ to $-e_n$ and preserves $e_j$ for $j \neq n$. 

Then the conjugation action of $\tW^!$ on $\tW$ sends integral elements to integral elements. 

\subsection{} We construct some polynomials $f_{w, C} \in \ZZ[v-v \i]$ (for $w$ an integral element in $\tW$ and $C$ a conjugacy class in $\tW$) as follows. 

If $w$ is a minimal element in the conjugacy class that contains it, then we set $f_{w, C}=\begin{cases} 1, & \text{ if } w \in C \\ 0, & \text{ if } w \notin C \end{cases}$. 

If $w$ is not a minimal element in the conjugacy class that contains it and that $f_{w', C}$ is already defined for all integral elements $w' \in \tW$ with $l(w')<l(w)$, then by \cite[Corollary 2.2]{He}, there exists $w_1 \approx w$ and $s \in \tS$ such that $l(s w_1 s)<l(w_1)=l(w)$, we define $f_{w, C}$ as $$f_{w, C}=(v-v \i) f_{w_1 s, C}+f_{s w_1 s, C}.$$

This completes the definition of $f_{w, C}$. One also sees from the definition that if $f_{w, C} \neq 0$, then all the coefficients of $f_{w, C}$ are nonnegative integer. 

By \ref{Hecke} (1) and \cite[Theorem 2.1]{He}, if $C$ is an integral conjugacy class in $\tW$ and $w, w' \in C_{\min}$, then $T_w \equiv T_{w'} \mod [H, H]$. Now we choose a representative $w_C \in C_{\min}$ for each integral conjugacy class $C$ in $\tW$. By \ref{Hecke} (1) \& (2) and \cite[Corollary 2.2]{He}, for any integral element $w \in \tW$, 

\[\tag{*} T_w \equiv \sum_C f_{w, C} T_{w_C} \mod [H, H].\] 

We call $f_{w, C}$ the class polynomials. 

Notice that the definition of $f_{w, C}$ depends  on the choice of the sequence of elements in $\tS$ used to conjugate $w$ to a minimal length element in its conjugacy class. We expect that $f_{w, C}$ is in fact, independent of such choice and is uniquely determined by the condition (*) above. This is true if one replaces $\tW$ by a finite Coxeter group and $H$ by the corresponding Hecke algebra (see \cite[Theorem 4.2]{GR}). 

\begin{cor}
Let $G$ be a classical simple algebraic group and $b \in G(L)$. Define $X^!_{w, \s}(b)=\{g I \in G(L)^!/I; g \i b \s(g) \in I \dot w I\}$. Let $C$ be an integral conjugacy class of $\tW$ and $w, w' \in C_{\min}$. Then $\dim(X^!_{w, \s_F}(b))=\dim(X^!_{w', \s_F}(b))$. 
\end{cor}

Proof. By \cite[5.1 (a)]{He}, there exists $x \in {}^S \tW$ and $v, v' \in W_{I(x)}$ such that $w \tilde \sim v x$ or $\iota(w) \tilde \sim v x$ and $w' \tilde \sim v' x$ or $\iota(w') \tilde \sim v' x$. Since $l(w)=l(w')$, we must have that $l(v)=l(v')$. By Corollary \ref{tt} and Lemma \ref{lang}, $\dim(X^!_{w, \s_F}(b))=\dim(X^!_{\iota(w), \s_F}(b))=\dim(X^^!_{v x, \s_F}(b))=\dim(X^!_{x, \s_F}(b))+l(v)$ and $\dim(X^!_{w', \s_F}(b))=\dim(X^!_{\iota(w'), \s_F}(b))=\dim(X^^!_{v' x, \s_F}(b))=\dim(X^!_{x, \s_F}(b))+l(v')$. Therefore $\dim(X^!_{w, \s_F}(b))=\dim(X^!_{w', \s_F}(b))$. The corollary is proved. \qed

\begin{thm}\label{class}
Let $G$ be a classical simple algebraic group and $b \in G(L)$. Let $w \in \tW$ be an integral element. Then 

$$\dim(X^!_{w, \s_F}(b))=\max_{C} \frac{1}{2}(l(w)-l(w_C)+\deg(f_{w, C}))+\dim(X^!_{w_C, \s_F}(b)),$$ where $C$ runs over integral conjugacy classes of $\tW$ and $w_C$ is a minimal length element in $C$.
\end{thm}

Proof. Let $C'$ be the integral conjugacy class that contains $w$. 

If $w \in C'_{\min}$, then $\frac{1}{2} \deg(f_{w, C})+\dim(X^!_{w_C, \s_F}(b) \neq -\infty$ if and only $f_{w, C} \neq 0$ and $X^!_{w_C, \s_F}(b) \neq \emptyset$, i.e., $C=C'$ and $b \in G(L) \cdot I \dot w_C I$. In this case, $f_{w, C}=1$ and $\deg(f_{w, C})=0$. By the previous Corollary, $\dim(X^!_{w, \s_F}(b))=\dim(X^!_{w_C, \s_F}(b))$. The theorem holds in this case. 

If $w \notin C'_{\min}$, we use the same sequence of elements in $\tS$ to conjugate $w$ to a minimal element in $C'$ as we did in the definition of $f_{w, C}$. Then there exists $w_1 \approx w$ and $s \in \tS$ such that $l(s w_1 s)<l(w_1)=l(w)$ and $f_{w, C}=(v-v \i) f_{w_1 s, C}+f_{s w_1 s, C}$. Hence $\deg(f_{w, C})=\max\{\deg(f_{w_1 s, C})+1, \deg(f_{s w_1 s, C})\}$ and $\frac{1}{2}(l(w)+\deg(f_{w, C}))=\max\{\frac{1}{2}(l(w_1 s)+\deg(f_{w_1 s, C}))+1, \frac{1}{2}(l(s w_1 s)+\deg(f_{s w_1 s, C}))+1\}$.

On the other hand, by Lemma \ref{indd}, $\dim(X^!_{w, \s_F}(b))=\dim(X^!_{w_1, \s_F}(b))=\max\{\dim(X^!_{w_1 s, \s_F}(b))+1, \dim(X^!_{s w_1 s, \s_F}(b))+1\}$. Now the theorem follows from induction on $l(w)$. \qed

\

By the same argument, one can prove the following result for $\s=\s_a$.

\begin{prop}\label{class2}
Let $G$ be a classical simple algebraic group and $b \in G(L)$. Let $w \in \tW$ be an integral element. If $X^!_{x, \s_a}(b)$ is finite dimensional for any $x \in \tW$ that is of minimal length in its conjugacy class in $\tW$, then $X^!_{w, \s_a}(b)$ is also finite dimensional. 
\end{prop}

%\begin{cor}
%upper bound of $\deg f_{w, C}$. 
%\end{cor}

\section{Some combinatorial properties}

\subsection{}  We follow the notations in \cite[1.4]{He}. 

A {\it double partition} $\tilde \l$ is a sequence $[(b_1, c_1), \cdots, (b_k, c_k)]$ with $(b_i, c_i) \in \NN \times \ZZ$ for all $i$ and $(b_1, c_1) \ge \cdots \ge (b_k, c_k)$ for the lexicographic order on $\NN \times \ZZ$. We write $|\tilde \l|$ for $(\sum b_i, \sum c_i)$ and $\l$ for $[b_1, \cdots, b_k]$. We set $|\O|=(0, 0)$ for the empty double partition $\O$. We call $\tilde \l$ {\it positive} if $c_i \ge 0$ for all $i$. We call $\tl$ {\it distinguished} if $b_i$ and $c_i$ are coprime for all $i$. We call $\tilde \l$ {\it special} if $c_i \in \{0, 1\}$ for all $i$. In this case, let $\underline{\tilde \l}$ be the double partition whose entries are $(b_1, 1-c_1), \cdots, (b_k, 1-c_k)$. Then $\underline {\tilde \l}$ is also special. 

Let $\cd \cp_A$ be the set of pairs of double partitions $(\tilde \l, \O)$ with $|\tilde \l|=(n, r)$ for some $0 \le r<n$. Let $\cd \cp$ be the set of pairs of double partitions $(\tilde l, \tilde \mu)$ such that $\tilde \mu$ is special and $|\tilde l|+|\tilde \mu|$ is of the form $(n, *)$ and $\cd \cp_{\ge 0}$ be the set of pairs $(\tilde \l, \tilde \mu) \in \cd \cp$ with $\tilde \l$ positive. 

Let $\tilde \l=[(b_1, c_1), \cdots, (b_k, c_k)]$ and $\tilde \mu=[(b_{k+1}, c_{k+1}), \cdots, (b_l, c_l)]$ with $(\tilde \l, \tilde \mu) \in \cd \cp$. We say that $(\tilde \l, \tilde \mu)$ is {\it distinguished} if $b_i$ and $c_i$ are coprime for all $i \le k$ and $c_i=1$ for $i>k$. 

If $\tW=\tW(A_{n-1})$, then any element in $\tW$ is integral. Then there is a bijection between the set of conjugacy classes of $\tW(A_{n-1})$ and $\cd \cp_A$. We denote the conjugacy class that corresponds to $(\tilde \l, \O)$ by $\co^A_{(\tilde \l, \O)}$. 

If $\tW=\tW(C_n)$, then an element in $\tW$ is integral if and only if it is in $W_a(C_n)$. Let $\sim$ be the equivalent relation on $\cd \cp$ defined by $(\tilde \l, \tilde \mu) \sim (\tilde \l, \underline{\tilde \mu})$ for all $(\tilde \l, \tilde \mu) \in \cd \cp$. Then there is a bijection between the set of $\tW(C_n)$-conjugacy classes in $W_a(C_n)$ and $\cd \cp_{\ge 0}/\sim$. We denote the conjugacy class that corresponds to $(\tilde \l, \tilde \mu)$ by $\co^C_{(\tilde \l, \tilde \mu)}$. Then $\co^C_{(\tilde \l, \tilde \mu)}=\co^C_{(\tilde \l, \underline{\tilde \mu})}$. 

\subsection{} Some representatives of the conjugacy classes are described as follows. 

Let $W(C_n)=(\ZZ/2 \ZZ)^n \rtimes S_n$ be the set of permutations $\s$ on $\{\pm 1, \cdots, \pm n\}$ with $\s(-i)=-\s(i)$ for all $i$. If $\s \in W(C_n)$ and there is only one or two orbits on $\{\pm 1, \cdots, \pm n\}$ consisting more than one element and the  orbit(s) are of the form $$i_1 \to i_2 \to \cdots \to i_k \to i_1 \text{ and/or } -i_1 \to -i_2 \to \cdots \to -i_k \to -i_1,$$ then we simply write $(i_1 i_2 \cdots i_k)$ for $\s$. For a pair of partitions $(\l, \mu)$ with $\l=[a_1, \cdots, a_k]$ and $\mu=[a_{k+1}, \cdots, a_l]$ and $\sum_{i=1}^l a_i=n$, we set  \begin{align*} w_{(\l, \mu)}=& (|\l|+a_{k+1}, -|\l|-a_{k+1}) (|\l|+a_{k+1}+a_{k+2}, -|\l|-a_{k+1}-a_{k+2}) \cdots (n, -n) \\ & (1 2 \cdots a_1) \cdots (n-a_l+1, n-a_l+2, \cdots, n). \end{align*}

Let $\tilde \l=[(b_1, c_1), \cdots, (b_k, c_k)]$ and $\tilde \mu=[(b_{k+1}, c_{k+1}, \cdots, (b_l, c_l)]$ with $(\tilde \l, \tilde \mu) \in \cd \cp$. Set $w_{(\tl, \tilde \mu)}^{st}=\e^{[a_1, \cdots, a_n]} w_{(\l, \mu)}$, where $w_{(\l, \mu)}$ is defined as above and $$a_{b_1+\cdots+b_{j-1}+i}=\lceil \frac{(b_j+1-i) c_j}{b_j} \rceil-\lceil \frac{(b_j-i) c_j}{b_j} \rceil$$ for $1 \le j \le l$ and $1 \le i \le b_j$. Here $\lceil x \rceil=\min\{n \in \ZZ; n \ge x\}$. 

If $(\tilde \l, \O) \in \cd \cp_A$, we call $w^{st}_{(\tilde \l, \O)}$ the standard representative of $\co^A_{(\tilde \l, \O)}$. If $(\tilde \l, \tilde \mu) \in \cd \cp_{\ge 0}$, we call $w^{st}_{(\tilde \l, \tilde \mu)}$ the standard representative of $\co^C_{(\tilde \l, \tilde \mu)}$. If moreover, $(\tilde \l, \tilde \mu)$ is distinguished, then by \cite[5.5]{He}, $S_n \cdot w^{st}_{(\tl, \tilde \mu)} \cap {}^S \tW$ (here $S_n$ acts on $\tW$ by conjugation) contains a unique element. We denote this element by $w^f_{(\tl, \tilde \mu)}$ and call it the {\it fundamental element} associated to $(\tl, \tilde \mu)$. 

\subsection{}\label{pw} Let $\tilde \l=[(b_1, c_1), \cdots, (b_k, c_k)]$ and $\tilde \mu=[(b_{k+1}, c_{k+1}), \cdots, (b_l, c_l)]$ with $(\tilde \l, \tilde \mu) \in \cd \cp$. Set $d_i=\gcd(b_i, c_i)$ and $m=\sum_{l'>k, c_{l'}=0} b_{l'}$. Let $\tilde \l'$ be the double partition whose entries are $(\frac{b_1}{d_1}, \frac{c_1}{d_1}), \cdots, (\frac{b_1}{d_1}, \frac{c_1}{d_1})$ (with $d_1$-times), $\cdots$, $(\frac{b_k}{d_k}, \frac{c_1}{d_k}), \cdots, (\frac{b_k}{d_k}, \frac{c_k}{d_k})$ (with $d_k$-times) and $(1, 0), (1, 0), \cdots, (1, 0)$ (with $m$-times) and $\tilde \mu'$ be the double partition whose entries are $(b_{l'}, c_{l'})$ with $l'>k$ and $c_{l'}=1$. Set $d(\tilde \l, \tilde \mu)=(\tilde \l', \tilde \mu')$. Then $d(\tilde \l, \tilde \mu)$ is distinguished and $w^f_{d(\tilde \l, \tilde \mu)}$ is defined. Set $d'(\tilde \l, \tilde \mu)=d(\tilde \l', \underline{\tilde \mu'})$. Then it is easy to see that  $d'(\tilde \l, \tilde \mu)$ is of the form $(*, \O)$. 

Recall the map $f: \tW \to B(\tW)$ defined in subsection \ref{ne}. If $\tW=\tW(A_{n-1})$, then $\tW=\sqcup_{(\tilde \l, \O) \in \cd \cp_A} \co^A_{(\tilde \l, \O)}$. Moreover, for any $(\tilde \l, \O) \in \cd \cp_A$, $w \in \co^A_{(\tilde \l, \O)}$ and $w' \in \co^A_{d(\tilde \l, \O)}$, we have that $f(w)=f(w')$. Assume that $\tilde \l=[(b_1, c_1), \cdots, (b_k, c_k)]$, then $[v_{w}]=[v_{w'}]$ is the $S_n$-orbit of $(e_1, \cdots, e_n)$, where $e_i=\frac{c_j}{b_j}$ for $b_1+\cdots+b_{j-1}<i \le b_1+\cdots+b_j$. Now for any distinguished double partition $\tl$ with $|\tl|=(n, r)$ for some $0 \le r \le n-1$, set $[\tl]_A=\sqcup_{d(\tilde \l', \O)=(\tilde \l, \O)} \co^A_{(\tilde \l', \O)}$. Then $[\tl]_A$ is a fiber of the map $f$ and we have that $$\tW(A_{n-1})=\sqcup [\tl]_A.$$

If $\tW=\tW(C_n)$, then $W_a=\sqcup_{[(\tilde \l, \tilde \mu)] \in \cd \cp_{\ge 0}/\sim} \co^C_{(\tilde \l, \tilde \mu)}$. Moreover, for any $(\tilde \l, \tilde \mu) \in \cd \cp_{\ge 0}$, $w \in \co^C_{(\tilde \l, \tilde \mu)}$, $w' \in \co^C_{d(\tilde \l, \tilde \mu)}$ and $w'' \in \co^C_{d'(\tilde \l, \tilde \mu)}$, we have that $f(w)=f(w')=f(w'')$. Assume that $\tilde \l=[(b_1, c_1), \cdots, (b_k, c_k)]$, then $[v_{w}]=[v_{w'}]$ is the $W(C_n)$-orbit of $(e_1, \cdots, e_n)$, where $e_i=\frac{c_j}{b_j}$ for $b_1+\cdots+b_{j-1}<i \le b_1+\cdots+b_j$ and $e_i=0$ for $i>b_1+\cdots+b_k$. Now for any distinguished positive double partition $\tl$, set $[\tl]_C=\sqcup_{d'(\tilde \l', \tilde \mu)=(\tilde \l, \O)} \co^C_{(\tilde \l', \tilde \mu)}$. Then $[\tl]_C$ is a fiber of the map $f$ and we have that $$W_a(C_n)=\sqcup [\tl]_C.$$

The following results are proved in \cite[Theorem 5.2]{He}. 

\begin{thm}\label{minA}
Let $\tl$ be a distinguished double partition with $|\tl|=(n, r)$ for some $0 \le r \le n-1$. Then $w \to w^f_{(\tilde \l, \O)}$ for any $w \in \co^A_{(\tilde \l, \O)}$. In particular, $w^f_{(\tilde \l, \O)}$ is a minimal length element in $\co^A_{(\tilde \l, \O)}$ and for any minimal length element $w$ in $\co^A_{(\tilde \l, \O)}$, $w \approx w^f_{(\tilde \l, \O)}$. If moreover, $(\tilde \l', \O) \in \cd \cp_A$ with $d(\tilde \l', \O)=(\tilde \l, \O)$, then for any $w \in \co^A_{(\tl', \O)}$, there exists $x \in W_{I(w^f_{(\tl, \O)})}$ such that $w \tilde \to x w^f_{(\tl, \O)}$.
\end{thm}

\begin{thm}\label{minC}
Let $(\tl, \tilde \mu) \in \cd \cp_{\ge 0}$ be distinguished. Then $w \to w^f_{(\tilde \l, \tilde \mu)}$ for any $w \in \co^C_{(\tilde \l, \tilde \mu)}$. In particular, $w^f_{(\tilde \l, \tilde \mu)}$ is a minimal length element in $\co^C_{(\tilde \l, \tilde \mu)}$ and for any minimal length element $w$ in $\co^C_{(\tilde \l, \tilde \mu)}$, $w \approx w^f_{(\tilde \l, \tilde \mu)}$. If moreover, $(\tilde \l', \tilde \mu') \in \cd \cp_{\ge 0}$ with $d(\tilde \l', \tilde \mu')=(\tilde \l, \tilde \mu)$, then for any $w \in \co^C_{(\tl', \tilde \mu')}$, there exists $x \in W_{I(w^f_{(\tl, \tilde \mu)})}$ such that $w \tilde \to x w^f_{(\tl, \tilde \mu)}$.
\end{thm}

\begin{cor}\label{xx}
(1) Let $\tl$ be a distinguished double partition with $|\tl|=(n, r)$ for some $0 \le r \le n-1$. If $w$ is a minimal length element in $[\tl]_A$, then $w \in \co^A_{(\tl, \O)}$. 

(2) Let $\tl$ be a positive distinguished double partition. If $w$ is a minimal length element in $[\tl]_C$, then $w \in \co^C_{(\tl, \O)}$. 
\end{cor}

Proof. Let $(\tilde \l', \tilde \mu)$ be a pair of double partitions that represents the conjugacy class of $w$. By Theorem \ref{minA} and \ref{minC}, there exists $x \in W_f$ such that $w \tilde \to x w^f_{d(\tilde \l', \tilde \mu)}$ and $l(x w^f_{d(\tilde \l', \tilde \mu)})=l(x)+l(x w^f_{d(\tilde \l', \tilde \mu)})$. Since $f(w)=f(w^f_{d(\tilde \l', \tilde \mu)})$, we must have that $x=1$ and $d(\tilde \l', \tilde \mu)=(\tilde \l', \tilde \mu)$. In other words, $(\tilde \l', \tilde \mu)$ is distinguished. 

So if $\tW=\tW(A_{n-1})$, then $\tl'$ is distinguished and $\tl'=\tl$. In this case $w \in \co^A_{(\tl, \O)}$. If $\tW=\tW(C_n)$, then $(\tilde \l', \tilde \mu)$ and $(\tilde \l', \underline{\tilde \mu})$ are distinguished. Therefore $\tilde \mu=\O$ and $\tl'=\tl$ is distinguished. In this case, $w \in \co^C_{(\tl, \O)}$. \qed

\begin{cor}\label{good2}
Let $w \in \tW(A_{n-1})$ (resp. $w \in W_a(C_n)$). Then $w$ is a good element if and only if $w$ is a minimal length element in $\co^A_{(\tl, \O)}$ (resp. $\co^C_{(\tl, \O)}$) for some distinguished double partition $\tl$. 
\end{cor}

Proof. This follows from the previous Corollary and Proposition \ref{good}. \qed

\begin{cor}\label{good3}
Let $\tW=\tW(A_{n-1})$ or $\tW(C_n)$. Let $w, w' \in \tW_{int}$ be good elements. Then $w$ and $w'$ are in the same fiber of the map $f: \tW_{int} \to B(\tW)$ if and only if $w \approx w'$. 
\end{cor}

Proof. This follows from the previous Corollary, Theorem \ref{minA} and Theorem \ref{minC}. \qed

\subsection{}\label{par} Let $\tW=\tW(A_{n-1})$ or $\tW(C_n)$. Set $B(\tW_{int})=f(\tW_{int})$. For $\aa \in B(\tW_{int})$ and $w \in \tW_{int}$, we write $\aa \preceq w$ if there exists $w' \in f \i(\aa)_{\min}$ with $w' \le w$ for the Bruhat order in $\tW$. By \cite[Lemma 4.4]{He072}, if $w_1 \tilde \to w$ and $\aa \preceq w$, then $\aa \preceq w_1$. For $\aa, \aa' \in B(\tW_{int})$, we write $\aa' \preceq \aa$ if there exists $w \in f\i(\aa)_{\min}$ with $\aa' \preceq w$. By \cite[Corollary 2.4]{He}, $\preceq$ is a partial order on $B(\tW_{int})$.

\section{The main result}

\begin{lem}\label{eq1}
Let $G=PGL_n$ or $PSP_{2 n}$. Let $\aa \in B(\tW_{int})$. Then 

(1) For any good elements $w, w' \in f \i(\aa)$, we have that $$G(L) \cdot_\s I \dot w' I=G(L) \cdot_\s I \dot w I=G(L)' \cdot_\s I \dot w I.$$  Now we define $X_{\aa, \s}=G(L) \cdot_\s I \dot w I$ for any good element $w \in \aa$. 

(2) Let $O \subset f \i(\aa)$ be a conjugacy class of $\tW$ and $x \in O_{\min}$. Then $G(L) \cdot_\s I \dot x I \subset X_{\aa, \s}$. 
\end{lem}

Proof. (1) By Corollary \ref{good3}, $w \approx w'$. Thus by Lemma \ref{red2}, $G(L) \cdot_\s I \dot w' I=G(L) \cdot_\s I \dot w I$. Notice that $G(L)=\sqcup_{\t \in \tW, l(\t)=0} G(L)' \dot \t$. Thus \begin{align*} G(L) \cdot_\s I \dot w I &=\cup_{\t \in \tW, l(\t)=0} G(L)' \cdot_\s \dot \t I \dot w I \s(\dot \t) \i \\ &=\cup_{\t \in \tW, l(\t)=0} G(L)' \cdot_\s I \dot \t w \dot \t \i I.\end{align*} Since $l(\t w \t \i)=l(w)$, $\t w \t \i$ is also a good in $f \i(\aa)$.  By Corollary \ref{good3} and Lemma \ref{red1}, $G(L)' \cdot_\s I \dot w I=G(L)' \cdot_\s I \dot \t \dot w \dot \t \i I$. Hence $G(L) \cdot_\s I \dot w I=G(L)' \cdot_\s I \dot w I$. 

(2) If $G=PGL_n$, then $f \i(\aa)=[\tilde \l]_A$ for some distinguished double partition $\tl$. Hence $O$ is of the form $\co^A_{\tl', \emptyset}$ with $d(\tl', \emptyset)=(\tl, \emptyset)$. By Theorem \ref{minA}, $x \tilde \approx u w^f_{(\tl, \O)}$ for some $u \in W_{I(w^f_{(\tl, \O)})}$. Now the statement follows from Lemma \ref{red2} and \ref{red3}. 

If $G=PGL_{2n}$, then $f \i(\aa)=[\tl]_C$ for some positive distinguished double partition $\tl$. Hence $O$ is of the form $\co^C_{(\tilde \l', \tilde \mu)}$ with $d'(\tl', \tilde \mu)=(\tl, \emptyset$. Set $d(\tilde \l', \tilde \mu)=(\tilde \l'', \tilde \mu')$. Then by Theorem \ref{minC}, $x \tilde \approx u w^f_{(\tl'', \tilde \mu')}$ for some $u \in W_{I(w^f_{(\tl'', \tilde \mu')})}$. By Lemma \ref{red2} and \ref{red3}, $G(L) \cdot_\s I \dot x I=G(L) \cdot_\s I \dot u \dot w^f_{(\tl'', \tilde \mu')} I \subset G(L) \cdot_\s I \dot w^f_{(\tl'', \tilde \mu')} I$. Since $w^f_{(\tl'', \tilde \mu')}$ is a minimal length element in $\co^C_{(\tilde \l'', \tilde \mu')}=\co^C_{(\tilde \l'', \underline{\tilde \mu'})}$, by the same argument, we have that $G(L) \cdot_\s I \dot w^f_{(\tl'', \tilde \mu')} I \subset G(L) \cdot_\s I \dot w^f_{d(\tl'', \underline{\tilde \mu'})} I=X_{\aa, \s}$. \qed 

\begin{prop}\label{par1}
Let $G=PGL_n$ or $PSP_{2 n}$ and $w \in \tW_{int}$. Then $$\overline{G(L) \cdot_\s I \dot w I}=G(L) \cdot_\s \overline{I \dot w I}=\sqcup_{\aa \preceq w} X_{\aa, \s}.$$
\end{prop}

Proof. By the proof of \cite[Prop 18]{V}, for any $w \in \tW$, $\overline{G(L) \cdot_\s I \dot w I}=G(L) \cdot_\s \overline{I \dot w I}$. By Lemma \ref{dg}, the union is in fact a disjoint union. If $\aa \preceq w$, then there exists a minimal length element $w'$ in $f \i(\aa)$ such that $w' \le w$. Hence $X_{\aa, \s}=G(L) \cdot_\s I \dot w' I \subset G(L) \cdot_\s \overline{I \dot w I}$. Now we prove that $G(L) \cdot_\s I \dot w I \subset \cup_{\aa \preceq w} X_{\aa, \s}$ by induction on $l(w)$. 

If $w$ is not a minimal length element in the conjugacy class of $\tW$ that contains it, then by Theorem \ref{minA} and \ref{minC} there exists $w_1 \approx w$ and $w_1 \xrightarrow i w_2$ with $l(w_2)<l(w_1)$. In particular, $w_2<w_1$ and $s_i w_1<w_1$. By induction hypothesis, \begin{align*} G(L) \cdot_\s I \dot w I &=G(L) \cdot_\s I \dot w_1 I=G(L) \cdot_\s I \dot w_2 I \cup G(L) \cdot_\s I \dot s_i \dot w_1 I  \\ & \subset \cup_{\aa \preceq w_1} X_{\aa, \s}.\end{align*} By subsection \ref{par}, $\aa \preceq w$ if and only if $\aa \preceq w_1$. So $G(L) \cdot_\s I \dot w I \subset \cup_{\aa \preceq w} X_{\aa, \s}$. 

If $w \in f \i(\aa)$ is a minimal length element in its conjugacy class, then by Theorem \ref{minA} and \ref{minC}, $\aa\preceq w$. By Lemma \ref{eq1}, $G(L) \cdot_\s I \dot w I \subset X_{\aa, \s}$. \qed

\begin{cor}
Let $G=PGL_n$ or $PSP_{2 n}$ and $\aa \in B(\tW_{int})$. Then $\overline{X_{\aa, \s}}=\sqcup_{\aa' \preceq \aa} X_{\aa', \s}$. 
\end{cor}

\begin{rmk}
If $\s=\s_F$ is a Frobenius morphism, then $X_{\aa, \s}$ is a single $\s$-conjugacy class and any $\s$-conjugacy class is of this form. In this case, the closure of $X_{\tl, \s}$ is a union of other $\s$-conjugacy classes and the explicit closure relation is obtained by Viehmann \cite[Prop 18]{V}. However, if $\s=\s_a$, then in general $X_{\aa, \s}$ contains infinitely many $\s$-conjugacy classes. 
\end{rmk}

\begin{prop}\label{str}
(1) Let $G=PGL_n$, then $G(L)=\sqcup_{\aa \in B(\tW)} X_{\aa, \s}$ is a stratification of $G(L)$. 

(2) Let $G=PSP_{2n}$, then $G(L)'=\sqcup_{\aa \in B(\tW_{int})} X_{\aa, \s}$ is a stratification of $G(L)'$. 
\end{prop}

\begin{rmk}
If $\s=\s_F$, then both parts follows from Kottwitz's classification of $\s$-conjugacy classes \cite{Ko97}. See subsection \ref{reK}. 
\end{rmk}

Proof. Notice that \begin{gather*} G(L)=\sqcup_{w \in \tW} I \dot w I=\cup_{w \in \tW} G(L) \cdot_\s I \dot w I=\cup_{w \in \tW} G(L) \cdot_\s \overline{I \dot w I}; \\ G(L)'=\sqcup_{w \in W_a} I \dot w I=\cup_{w \in W_a} G(L) \cdot_\s I \dot w I=\cup_{w \in W_a} G(L) \cdot_\s \overline{I \dot w I}.\end{gather*} Now the proposition follows from Proposition \ref{par1} and the previous Corollary.  \qed

\begin{prop}\label{class3}
Let $G=PGL_n$ or $PSP_{2n}$. Let $\aa \in B(\tW_{int})$. Then for any $w \in \tW$ and $b \in X_{\aa, \s}$, we have that 

(1) $\dim (X_{w, \s}(b))=\max_{C} \frac{1}{2}(l(w)+l(C)+\deg(f_{w, C}))-l(f \i(\aa))$ if $\s=\s_F$, here $C$ runs over conjugacy class of $\tW$ in $f \i(\aa)$.

(2) $\dim(X_{w, \s}(b))< \infty$ if $\s=\s_a$ with $a$ not a root of unity.
\end{prop}

Proof. It is easy to see that if $X_{w, \s}(b) \neq \emptyset$, then $w \in \tW_{int}$. Now let $C$ be an integral conjugacy class of $\tW$. If $X_{w_C, \s}(b) \neq \emptyset$, then $b \in G(L) \cdot_\s I \dot w_C I$. By Lemma \ref{eq1} and Lemma \ref{dg}, we must have that $C \subset f \i(\aa)$. 

Case I: $G=PGL_n$. Then $f \i(\aa)=[\tilde \l]_A$ for some distinguished double partition $\tl$. Hence $O$ is of the form $\co^A_{(\tl', \emptyset)}$ with $d(\tl', \emptyset)=(\tl, \emptyset)$. By Theorem \ref{minA}, $w_C \tilde \approx u w^f_{(\tl, \O)}$ for some $u \in W_{I(w^f_{(\tl, \O)})}$. By Lemma \ref{indd}, $\dim(X_{w_C, \s}(b))=\dim(X_{u w^f_{(\tl, \O)}, \s}(b))$. 

If $\s=\s_F$, then by Lemma \ref{lang} and Prop \ref{00}, $\dim(X_{u w^f_{(\tl, \O)}, \s}(b))=\dim(X_{w^f_{(\tl, \O)}, \s}(b))+l(u)=l(u)=l(w_C)-l(w^f_{(\tl, \O)})=l(C)-l(f \i(\aa))$. Now by Theorem \ref{class}, $\dim (X_{w, \s}(b))=\max_{C} \frac{1}{2}(l(w)+l(C)+\deg(f_{w, C}))-l(f \i(\aa))$. 

If $\s=\s_a$ for some $a$ not a root of unity, then by Lemma \ref{lang} and Prop \ref{00}, $\dim(X_{u w^f_{(\tl, \O)}, \s}(b)) \le \dim(X_{w^f_{(\tl, \O)}, \s}(b)+l(w_0^S)<\infty$. Hence by Prop \ref{class2}, $\dim(X_{w, \s}(b))<\infty$. 

Case II: $G=PGL_{2n}$. Then $f \i(\aa)=[\tl]_C$ for some positive distinguished double partition $\tl$. Hence $O$ is of the form $\co^C_{(\tilde \l', \tilde \mu)}$ with $d'(\tl', \tilde \mu)=(\tl, \emptyset)$. Let $d(\tilde \l', \tilde \mu)=(\tilde \l'', \tilde \mu')$. Then by Theorem \ref{minC}, $x \tilde \approx u w^f_{(\tl'', \tilde \mu')}$ for some $u \in W_{I(w^f_{(\tl'', \tilde \mu')})}$. By Lemma \ref{indd}, $\dim(X_{w_C, \s}(b))=\dim(X_{u w^f_{(\tl'', \tilde \mu')}, \s}(b))$. 

If $\s=\s_F$, then by Lemma \ref{lang} and Prop \ref{00}, $\dim(X_{u w^f_{(\tl'', \tilde \mu')}, \s}(b))=\dim(X_{w^f_{(\tl'', \tilde \mu')}, \s}(b))+l(u)$.  Since $w^f_{(\tl'', \tilde \mu')}$ is a minimal length element in $\co^C_{(\tilde \l'', \tilde \mu')}=\co^C_{(\tilde \l'', \underline{\tilde \mu'})}$, by the same argument, we have that $\dim(X_{w^f_{(\tl'', \tilde \mu')}, \s}(b))=l(w^f_{(\tl'', \tilde \mu')})-l(f \i(\aa))$. Therefore $\dim(X_{w_C, \s}(b))=l(C)-l(f \i(\aa))$. Now by Theorem \ref{class}, $\dim (X_{w, \s}(b))=\max_{C} \frac{1}{2}(l(w)+l(C)+\deg(f_{w, C}))-l(f \i(\aa))$. 

If $\s=\s_a$ for some $a$ not a root of unity, then by Lemma \ref{lang} and Prop \ref{00}, $\dim(X_{w^f_{(\tl'', \tilde \mu')}, \s}(b)) \le \dim(X_{w^f_{(\tl, \O)}, \s}(b))+l(w_0^S)<\infty$. Again by Lemma \ref{lang} and Prop \ref{00}, $\dim(X_{u w^f_{(\tl'', \tilde \mu')}, \s}(b))<\infty$. Hence by Prop \ref{class2}, $\dim(X_{w, \s}(b))<\infty$. \qed

\subsection{} Let us come to the case where $\s=\s_1$ is the identity map. In this case, the $\s$-conjugacy classes in $G(L)$ are just the usual conjugacy classes. Let $b$ be a regular semisimple, integral element in $G(L)$. Here integral means the elements in $G(L) \cdot I$. It is shown by Kazhdan and Lusztig in \cite{KL} that $X_{(1, id)}(b)$ is finite dimensional and a conjectural dimension formula is also given there. If moreover, $b$ is elliptic (i.e., its centralizer is an anisotropic torus), then $X_{1, id}(b)$ has only finitely many irreducible components and is an algebraic variety. The conjecture is proved later by Bezrukavnikov in \cite{Be}. (Actually they considered only topologically unipotent elements, but the general can be reduced to that case using Jordan decomposition). Now by the same argument as we did above, one can show the following result. This answers the question in \cite[Section 7]{L3} for $G=PGL_n$ or $PSP_{2n}$. 

\begin{prop}
Let $G=PGL_n$ or $PSP_{2n}$. Let $b$ be a regular semisimple integral element in $G(L)$. Then for any $w \in \tW$, $X_{w, id}(b)$ is finite dimensional. If moreover, $b$ is elliptic, then $X_{w, id}(b)$ is an algebraic variety. 
\end{prop}

\

Now we can prove our main theorem. 

\begin{thm}\label{main}
Let $G=PGL_n$ or $PSP_{2n}$, $\aa \in B(\tW_{int})$ and $b \in X_{\aa, \s}$. If either (a) $\s=\s_F$ or (b) $\s=\s_a$ with $a$ not a root of unity or (c) $\s=id$ and $b$ is a regular semisimple integral element, then any simple $Z_{G(L), \s}(b)$-module occurs as a subquotient of $H^*_c(X_{w, \s}(b))$ for some $w \in \tW$ must occurs as a subquotient of $H^*_c(X_{x, \s}(b))$ for some minimal length element $x$ in an integral conjugacy class of $\tW$ on $f \i(\aa)$. More precisely, 

(1) If $G=PGL_n$ and $f \i(\aa)=[\tl]_A$, then we may take $x$ to be an element of the form $u w^f_{(\tl, \O)}$, where $u \in W_{I(w^f_{(\tl, \O)})}$. 

(2) If $G=PSP_{2n}$ and $f \i(\aa)=[\tl]_C$, then we may take $x$ to be an element of the form $u w^f_{(\tl', \tilde \mu')}$, where$(\tl', \tilde \mu')$ is distinguished with $d(\tl', \tilde \mu')=(\tl, \O)$ and $u \in W_{I(w^f_{(\tl', \tilde \mu')})}$. 
\end{thm}

Proof. By Lemma \ref{indd}, simple $Z_{G(L), \s}(b)$-module occurs as a subquotient of $H^*_c(X_{w, \s}(b))$ for some $w \in \tW$ must occurs as a subquotient of $H^*_c(X_{x, \s}(b))$ for some minimal length element $x$ in an integral conjugacy class of $\tW$. Now let $C$ be an integral conjugacy class and $x \in C_{\min}$. If $X_{x, \s}(b) \neq \emptyset$, then by Lemma \ref{eq1} and Lemma \ref{dg}, we must have that $C \subset f \i(\aa)$. The ``more precise'' part follows from Theorem \ref{minA} and Theorem \ref{minC}. \qed

\begin{cor}\label{finiteA}
Assume that $G=PGL_n$ is defined and split over $\FF_q$ and $F$ is the Frobenius morphism. Then any simple $G(\FF_q((\e)))$-module occurs as a subquotient of $H^*_c(X_{w, \s_F}(1))$ for some $w \in \tW$ must occurs as a subquotient of $H^*_c(X_{x, \s_F}(1))$ for some $x \in W_f$. 
\end{cor}

\begin{cor}\label{unique}
We keep the assumption in Theorem \ref{main}. Assume furthermore that $f \i(\aa)$ is a single conjugacy class, then any simple $Z_{G(L), \s}(b)$-module occurs as a subquotient of $H^*_c(X_{w, \s}(b))$ for some $w \in \tW$ must occurs as a subquotient of $H^*_c(X_{w^f_{(\tl, \O)}, \s}(b))$.
\end{cor}

%\begin{prop}
%(1) Let $G=PGL_n$ and $\tl$ be a distinguished double partition with $|\tl|=(n, r)$ for some $0 \le r \le n-1$. Then any affine Deligne-Lusztig representation (resp. unipotent affine Deligne-Lusztig representation) for $w^f_{(\tl, \O)}$ occurs as a simple subquotient of some $H^i_c(X_{w}(b), \cf_\th)$ (resp. $H^i_c(X_{w}(b), \bbq)$), where $w$ is a minimal length element in $\co^A_{(\tl', \O)}$ with $d(\tl', \O)=(\tl, \O)$. 

%(2) Let $G=PSP_{2n}$ and $\tl$ be a positive distinguished double partition. Then any affine Deligne-Lusztig representation (resp. unipotent affine Deligne-Lusztig representation) for $w^f_{(\tl, \O)}$ occurs as a simple subquotient of some $H^i_c(X_{w}(b), \cf_\th)$ (resp. $H^i_c(X_{w}(b), \bbq)$), where $w$ is a minimal length element in $\co^C_{(\tl', \tilde \mu)}$ with $d'(\tl', \tilde \mu)=(\tl, \O)$.
%\end{prop}

\subsection{} In the rest of this paper, we discuss in more details the special cases that $f \i(\aa)$ is a single conjugacy class. 

Let $\tl=[(b_1, c_1), \cdots, (b_k, c_k)]$ be a double partition. We call $\tl$ {\it super-distinguished} if $\tl$ is distinguished and all the entries $(b_1, c_1), \cdots, (b_k, c_k)$ are distinct. In this case, for $\tW=\tW(A_{n-1})$, $[\tl]_A=\co^A_{(\tl, \O)}$ and $I(w^f_{(\tl, \O)})=\emptyset$. 

We call $\tl$ {\it strictly positive} if $c_i>0$ for all $i$. Then for $\tW=\tW(C_n)$ and any super-distinguished strictly positive double partition $\tl$, $[\tl]_C=\co^C_{(\tl, \O)}$ and $I(w^f_{(\tl, \O)})=\emptyset$. 

Then it is easy to see that 

(a) For $G=PGL_n$, $f \i(\aa)$ contains a single conjugacy class if and only if $\aa=[\tl]_A$ for some super-distinguished double partition $\tl$.

(b) For $G=PSP_{2n}$, $f \i(\aa)$ contains a single conjugacy class if and only if $\aa=[\tl]_A$ for some super-distinguished strictly positive double partition $\tl$. 

\subsection{} Let $G=PGL_n$. Then for any dominant regular coweight $\chi$, $\e^\chi$ is of the form $w^f_{(\tl, \O)}$ for some super-distinguished double partition $\tl=[(1, a_1), \cdots, (1, a_n)]$, where $a_i-a_{i+1}=<\chi, \a_i> >0$. 

Similarly, let $G=PSP_{2n}$. Then for any dominant regular coweight $\chi$ that lies in the coroot lattice, $\e^\chi$ is of the form $w^f_{(\tl, \O)}$ for some super-distinguished strictly positive double partition $\tl=[(1, a_1), \cdots, (1, a_n)]$, where $a_i-a_{i+1}=<\chi, \a_i> >0$ and $a_n=<\chi, \a_n> >0$. 

\begin{cor}\label{ddd}
Let $\s=\s_F$ or $\s_a$ for $a$ not a root of unity. We assume that either 

(a) $G=PGL_n$ and $\chi$ is a dominant regular coweight or 

(b) $G=PSP_{2n}$ and $\chi$ is a dominant regular coweight that lies in the coroot lattice. 

Let $w \in \tW$. Then any simple representation of $Z_{G(L), \s}(\e^\chi)$ as a subquotient of $H^*_c(X_{w, \s}(\e^\chi))$ factors through $Z_{G(L), \s}(\e^\chi)/(Z_{G(L), \s}(\e^\chi) \cap I)$.
\end{cor}

\begin{rmk}
The special case for $\s=\s_F$ and $G=PGL_2$ or $PGL_3$ is obtained by Zbarsky \cite{Zb} in a different way. 
\end{rmk}

Proof. By Corollary \ref{unique}, any simple representation of $Z_{G(L), \s}(\e^\chi)$ occurs as a subquotient of $H^*_c(X_{w, \s}(\e^\chi))$ must occurs as a subquotient of $H^*_c(X_{\e^\chi, \s}(\e^\chi))$. By Prop \ref{00}, $X_{\e^\chi, \s}(\e^\chi) \subset T(L) I/I \cong T(L)/(T(L) \cap I)$. We also see from the proof of Prop \ref{00} that $Z_{G(L), \s}(\e^\chi) \subset T(L)$. Since $T(L) \cap I$ acts trivially on $T(L)/(T(L) \cap I)$, $Z_{G(L), \s}(\e^\chi) \cap I$ acts trivially on $X_{\e^\chi, \s}(\e^\chi)$ and also trivially on $H^*_c(X_{\e^\chi, \s}(\e^\chi))$. \qed

\subsection{} Now we consider another case. Let $G=PGL_n$ and $\t=w^f_{(\tl, \O)} \in \tW$, where $\tl=(n, r)$ for some $0<r<n$ with $r$ and $n$ are coprime. By \cite[Lemma 4.5]{He}, $\t=t^{\o_r} w_0^{{S-\{r\}}} w_0^S$. We call $\t$ a superbasic element. %In particular, $\t$ is a superbasic element and any superbasic element in $G(L)$ is $F$-conjugate to such a $\t$. 

\begin{lem}
We keep the notation as in the previous subsection. If $g \in G(L)$ such that $g \i \dot \t \s(g) \in I \dot \t$, then $g \in I$. 
\end{lem}

Proof. Assume that $g \in I \dot w I$ for $w \in \tW$. Then $g \in \dot \t \s(g) \dot \t \i I=I \dot \t \dot w \dot \t \i I$ and $I \dot w I \cap I \dot \t \i \dot w \dot \t I \neq \emptyset$. Therefore $w=\t w \t \i$. If $w \neq 1$, then there exist $i \in \tS$ such that $s_i w<w$. Since conjugation by $\t$ preserve the Bruhat order, we have that $s_{\t^n(i)} \le \t^n w \t^{-n}=w$. However, $\tS$ is a single $\t$-orbit. Therefore $s_j w<w$ for any $j \in \tS$. That is a contradiction. \qed

\

By the same method as we did in the proof of Corollary \ref{ddd}, we have that 

\begin{cor}\label{superbasic}
Let $G=PGL_n$ and $\t \in \tW$ a superbasic element. Let $w \in \tW$. Then any simple representation of $Z_{G(L), \s}(\dot \t)$ occurs as a subquotient of $H^*_c(X_{w, \s}(\dot \t))$ is trivial. 
\end{cor}

\section*{Acknowledgement} 
We thank Victor Ginzburg, Ulrich G\"ortz, George Lusztig, Michael Rapoport, Eva Viehmann and Yong-Chang Zhu for some helpful discussions. Part of this work was done during my visit at University of Bonn. We thank them for the warm hospitality. 

%\bibliography{heref}
\bibliographystyle{amsalpha}

\end{document}